\documentclass[11pt]{amsart}

\usepackage{amsfonts, enumerate} 
\usepackage{amssymb}
\usepackage[final]{graphicx, graphics,stmaryrd} 
\usepackage[margin=3cm]{geometry} 
\usepackage{amscd}
\usepackage[breaklinks,bookmarksopen,bookmarksnumbered]{hyperref}

\usepackage{tensor}
\usepackage[english]{babel}
\usepackage[utf8]{inputenc}
\usepackage{amsthm}
\usepackage{graphics}
\usepackage{amsmath}
\usepackage{amstext}
\usepackage{engrec}
\usepackage{rotating}
\usepackage[safe,extra]{tipa}
\usepackage{multirow}
\usepackage{enumitem}
\usepackage{bm}
\usepackage{stackrel} 
\usepackage{tikz-cd} 
\numberwithin{equation}{section}

\usepackage{pgf,tikz,pgfplots}
\usepackage{mathrsfs}
\usetikzlibrary{arrows}

\newtheorem{teo}{Theorem}[subsection]

\newtheorem{prop}[teo]{Proposition}
\newtheorem{lemma}[teo]{Lemma}

\newtheorem*{teon}{Theorem}

\theoremstyle{definition}
\newtheorem{defin}[teo]{Definition}
\newtheorem{rmk}[teo]{Remark}

\newcommand{\mb}{\mathbb}
\newcommand{\mc}{\mathcal}
\newcommand{\msf}{\mathsf}
\newcommand{\mfk}{\mathfrak}

\def\Oo{\mathcal O}
\def\R{\mathcal{R}}
\def\C{\mathbb C}

\def\Z{\mathbb Z}

\def\Q{\mathbb Q}
\def\p{\mathbb P}
\def\cc{\mathcal C}
\def\E{\mathbb E}
\def\M{\mathcal{M}}
\def\A{\mathbb A}

\def\mm{\bm{\mathsf{M}}}
\def\rr{\bm{\mathsf{R}}}
\def\U{\mathcal U}
\def\cc{\mathcal{C}}
\def\D{\mathcal{D}}

\def\I{\mc{I}}

\newcommand{\ssC}{\mathsf{C}}

\newcommand{\m}{\mbox}

\newcommand{\cor}{\textit}

\newcommand{\fine}{\qed\newline}
\newcommand{\xx}{\otimes}

\newcommand{\beginCD}{\begin{equation*}\begin{CD}}
\newcommand{\enCD}{\end{CD}\end{equation*}}

\newcommand{\G}{\mfk G}

\newcommand{\we}{\wedge}

\DeclareMathOperator{\age}{age}
\DeclareMathOperator{\fix}{Fix}
\DeclareMathOperator{\defo}{Def}
\DeclareMathOperator{\Aut}{Aut}
\DeclareMathOperator{\diag}{Diag}

\DeclareMathOperator{\im}{Im}

\DeclareMathOperator{\sing}{Sing}

\DeclareMathOperator{\Ho}{Hom}

\DeclareMathOperator{\Spec}{Spec}

\DeclareMathOperator{\GL}{GL}

\DeclareMathOperator{\QR}{QR}
\DeclareMathOperator{\Pic}{Pic}

\DeclareMathOperator{\rk}{rk}
\DeclareMathOperator{\nc}{nc}

\DeclareMathOperator{\Bl}{Bl}

\def\rr{\bm{\mathsf{R}}}

\def\lq{\llbracket}
\def\rq{\rrbracket}

\newcommand{\nor}{\msf {nor}}

\DeclareMathOperator{\reg}{reg}

\DeclareMathOperator{\bal}{bal}

\DeclareMathOperator{\adm}{Adm}
\def\E{\mb E}
\def\I{\mb I}
\def\rr{\bm{\mathsf{R}}}
\def\mmu{\bm{\mu}}
\DeclareMathOperator{\inv}{inv}

\DeclareMathOperator{\gen}{gen}

\DeclareMathOperator{\Def}{Def}
\DeclareMathOperator{\Fun}{Fun}
\DeclareMathOperator{\Tor}{Tor}
\DeclareMathOperator{\Sym}{Sym}
\DeclareMathOperator{\ev}{ev}
\def\t{\tau}
\DeclareMathOperator{\ch}{ch}

\DeclareMathOperator{\sub}{Sub}

\providecommand{\Keywords}[1]{\textbf{Keywords: }#1}

\title{Birational geometry of moduli of curves with an $S_3$-cover}
\author{Mattia Galeotti}
\address{Mattia Galeotti,\ Università di Bologna,\ Piazza di Porta S. Donato 5, 40126 Bologna, Italy}
\email{galeotti.mattia.work@gmail.com}

\begin{document}

\maketitle

\begin{abstract}
We consider the space $\R_{g,S_3}^{S_3}$ of curves with a connected $S_3$-cover,
proving that for any odd genus $g\geq 13$ this moduli is of general type. Furthermore
we develop a set of tools that are essential in approaching the case of $G$-covers for any finite group $G$.
\end{abstract}


\section{Introduction}
The goal of this paper, alongside its companion \cite{gale19}, is to analyze the birational geometry
of the moduli space of curves equipped with a $G$-cover, where $G$ is any finite group.
More specifically, here we prove that the moduli space $\overline\R_{g,S_3}^{S_3}$ of genus $g$ connected twisted $S_3$-covers
is of general type for any odd genus $g\geq 13$.

In a series of seminal works, Eisenbud, Harris and Mumford proved that
the moduli $\M_g$ of genus $g$ smooth curves is a variety of general type
for any genus $g>23$ (see \cite{harmum82, har84, eihar87}).
The behavior is different for low genus, and in fact there is a large
literature showing that $\M_g$ is unirational for $g\leq 14$ and rationally
connected for $g=15$ (see for instance \cite{arbacor81, ver05, bruver05}).
Recently, Farkas, Jensen and Payne
showed 
that also $\M_{22}$ and $\M_{23}$ are of generaly type \cite{fjp20}.

Many finite covers of $\M_g$ with a modular interpretation, have the same
behavior: there is a transition genus such that
for any higher value the variety is of general type. The interesting
property is that in many cases the
transition genus is strictly lower. Farkas and Ludwig proved in \cite{farlud10} that 
the space $\R_{g,\mmu_2}^{\mmu_2}$ of curves with a non-trivial $2$-torsion bundle
is of general type for any genus $g\geq14$ and $g\neq 15$; Chiodo, Eisenbud, Farkas and Schreyer (see~\cite{cefs13})
proved the same in the case of $3$-torsion bundles for $g\geq 12$;
Farkas and Verra approached in \cite{farver10} the case of odd spin curves,
which is of general type for $g\geq12$.\newline

In this paper we focus on the moduli $\R_{g,G}$ of curves
with a principal $G$-bundle for $G$ finite group. We build its compactification $\overline\R_{g,G}$
as the moduli of twisted $G$-covers, that is balanced  representable maps
$\phi\colon\ssC\to BG$ where
$\ssC$ is a Deligne-Mumford stack such that its coarse space is a stable curve $C$,
with non-trivial stabilizer only at some nodes.
We properly introduce twisted $G$-covers in \S\ref{int_twist},
and in \S\ref{int_adm} we recall the equivalence
with the notion of admissible $G$-cover. 
In order to evaluate the Kodaira dimension of $\overline \R_{g,G}$
it is necessary to consider a desingularization $\widehat \R_{g,G}\to \overline\R_{g,G}$,
but it is difficult to evaluate the dimension of the space of $n$-canonical sections
$H^0(\widehat \R_{g,G},nK)$. Following the approach of \cite{harmum82},
we intend to focus in those cases where it is the same to evaluate the dimension of the space of $n$-canonical
sections over the regular locus of $\overline \R_{g,G}$, that is those cases such that
$$H^0\left(\widehat\R_{g,G},nK_{\widehat\R_{g,G}}\right)=H^0\left(\overline\R_{g,G}^{\reg},nK_{\overline\R_{g,G}^{\reg}}\right).$$

A pluricanonical form defined locally over $\overline \R^{\reg}_{g,G}$, extends uniquely
 to $\widehat \R_{g,G}$ in smooth points and in the so called canonical singularities (for
 a treatment of canonical and non-canonical singularities see \cite{reid80}).
Therefore, we need an analysis of the locus of non-canonical singularities,
as done by Chiodo and Farkas in \cite{chiofar12} for
level $\ell$ curves (that is curves with an $\ell$-torsion bundle), and by the author in \cite{gale15}
for the moduli space $\R_{g,\ell}^k$ with a line bundle $L$ such that $L^{\xx\ell}\cong \omega^{\xx k}$.

We search for a result allowing a global unique extension of pluricanonical forms.
Consider the natural project $\pi\colon \overline \R_{g,G}\to \overline\M_g$,
we know from the author work \cite{gale19} that
the locus of non-canonical singularities is bipartitioned as $\sing^{\nc}\overline\R_{g,G}=T_{g,G}\cup J_{g,G}$,
where the $T$-locus is defined as
$T_{g,G}:=\sing^{\nc}\overline\R_{g,G}\cap\pi^{-1}\sing^{\nc}\overline\M_g$,
while $J_{g,G}$ is the locus of ``new'' singularities that are not in the preimage
of $\overline\M_g$ singularities. We recall that the $J$-locus is empty for~$\overline\R_{g,S_3}$, a fact
proven in \cite[Theorem 5.2.1]{gale19}. 

In Theorems \ref{teoh0} and \ref{teo_extpcf}
we prove the following.

\begin{teon}
Consider a desingularization $\widehat \R_{g,G}\to\overline\R_{g,G}$.
If $G$ is abelian and $J_{g,G}=\varnothing$,
then
$$H^0\left(\widehat\R_{g,G},nK_{\widehat\R_{g,G}}\right)=H^0\left(\overline\R_{g,G}^{\reg},nK_{\overline\R_{g,G}^{\reg}}\right)$$
for $n$ sufficiently big and divisible. The same is true when $G=S_3$ is the symmetric group of order $3$.
\end{teon}

The extension result is obtained via a non-trivial generalization of Harris-Mumford techniques
used in \cite{harmum82}, which is particularly tricky in the case $G=S_3$ because the covering
$T_{g,S_3}\to \sing^{\nc}\overline\M_g$ is non-étale.\newline

Finally, the main result of this paper concerns the component of $\overline\R_{g,S_3}$ parametrizing connected
twisted $S_3$-covers, that we denote by $\overline\R_{g,S_3}^{S_3}$. The other
 components are finite covers of the well known moduli of level curves of order $2$ and $3$.
In order to prove the bigness of the canonical divisor of $\overline \R_{g,S_3}^{S_3}$, for every odd genus $g=2i+1$ we write it down as 
a linear combination
$$K_{\overline\R_{g,S_3}^{S_3}}=\alpha\cdot \overline\U_g+\beta\cdot (\pi^*\M^1_{g,i+1})+E +\gamma\cdot \lambda\in\Pic_{\Q}(\overline\R_{g,S_3}^{S_3}),$$
where $\alpha,\beta,\gamma$ are real coefficients, $E$ is a boundary effective divisor and $\lambda$ is the Hodge class,
which is known to be big.
Furthermore, the divisor $\U_g$ is the jumping locus for the Koszul cohomology $K_{i,1}$
of a particular rank $2$ vector bundle (see Definition \ref{def_U}), its effectiveness is proved in Proposition \ref{propeff},
while $\pi^*\M^1_{g,i+1}$ is the lifting of an effective $\M_g$ divisor introduced
by Harris and Mumford (see Definition \ref{def_mmm}). Therefore $\overline \U_g$
and $\pi^*\M^1_{g,i+1}$ are both effective.
In Proposition \ref{prop_fin} we prove that 
the coefficients $\alpha,\beta,\gamma$ can be chosen all positive for $i>5$, 
and this implies the bigness of $\lambda$ and the space $\overline\R_{g,S_3}^{S_3}$ being
of general type.
\begin{teon}
The moduli space $\overline\R_{g,S_3}^{S_3}$ is of general type for every odd genus $g\geq 13$.\newline
\end{teon}

In \S\ref{s2} we introduce the compactification of $\R_{g,G}$ via twisted $G$-covers and admissible $G$-covers.
In \S\ref{struct} we describe the boundary of $\overline\R_{g,S_3}$ and evaluate its canonical divisor.
In \S\ref{ext} the extension of pluricanonical forms is shown and in \S\ref{finally} we build the Koszul divisor
and conclude the proof.

\section*{Acknowledgements}
This work is the completion of my PhD thesis at \cor{Institut de Mathématiques de Jussieu}, thus  I want to thank
my PhD supervisor Alessandro Chiodo, for his patience and his advices. Furthermore, I want to thank
Michele Bolognesi and Gavril Farkas for the careful reading and the comments on my thesis work. Finally, I am grateful
to Claudio Fontanari, my post-doc tutor at \cor{Università degli Studi di Trento}, for his important suggestions,
and to Roberto Pignatelli for the fruitful conversations.

\section{Twisted and admissible covers}\label{s2}
In this section we recall the notions of twisted $G$-cover and admissible $G$-cover,
and their equivalence in Theorem \ref{teo_3e}.
In the following of the manuscript, we will adopt preferably the twisted cover notation in Section \ref{ext} and \ref{finally},
and the admissible cover notation in Section \ref{struct}.

In \S\ref{int_adm} we also
 introduce the notion of admissible $\mc H$-cover, where $\mc H$
is any subgroup conjugacy class of $G$. This concept is central in
classifying the components and boundary divisors of any moduli space $\overline\R_{g,G}$,
as we do in Section \ref{struct}, in particular for the case $G=S_3$.
In order to better understand these subjects, we also recall a series of group action tools.

Finally, in \S\ref{sec_monodromy} we introduce a description of admissible $G$-covers via a monodromy-type
approach. This description allows to treat admissible $G$-covers as particular classes of
group morphisms.

\subsection{Introduction to twisted $G$-covers}\label{int_twist}

\begin{defin}[Twisted curve]
A twisted $n$-marked $S$-curve is a diagram 
$$\begin{array}{ccc}
\Sigma_1,\Sigma_2,\dots,\Sigma_n & \subset &\ssC \\
 & & \downarrow \\
 & & C\\
 & & \downarrow \\
 & & S.
 \end{array}$$
 \begin{enumerate}
 \item $\ssC$ is a Deligne-Mumford stack, proper over $S$,
 and étale locally it is a nodal curve over~$S$;
 \item the $\Sigma_i\subset \ssC$ are disjoint closed substacks in the smooth
 locus of $\ssC\to S$ for all~$i$;
 \item $\Sigma_i\to S$ is an étale gerbe for all~$i$;
 \item $\ssC\to C$ exhibits $C$ as the coarse space of $\ssC$, and it is an isomorphism over $C_{\gen}$.
 \end{enumerate}
\end{defin}

We recall that, given a scheme $U$ and a finite abelian group $\mmu$ acting on $U$,
the stack $[U\slash \mmu]$ is the category of principal $\mmu$-bundles $E\to T$, for
any scheme $T$, equipped with a $\mmu$-equivariant morphism $f\colon E\to U$.
The stack $[U\slash \mmu]$ is a proper Deligne-Mumford stack and has a natural
morphism to its coarse scheme $U\slash\mmu$.\newline

By the definition of twisted curve we get the local pictures:
\begin{itemize}
\item \cor{At a marking}, morphism $\ssC\to  C\to S$ is locally isomorphic to 
$$\left[\Spec A[x']\slash \mmu_r\right]\to \Spec A[x]\to\Spec A$$
for some normal ring $A$ and some integer $r>0$.
Here $x=(x')^{r}$, and
$\mmu_r$ is the cyclic group of order $r$ acting on $\Spec A[x']$ by the action 
$\xi\colon x'\mapsto \xi x'$ for any $\xi\in\mmu_r$.\newline

\item \cor{At a node}, morphism $\ssC \to C\to  S$ is locally isomorphic to 
$$\left[\Spec\left( \frac{A[x',y']}{(x'y'-a)}\right)\slash \mmu_r\right]\to \Spec\left( \frac{A[x,y]}{(xy-a^\ell)}\right)\to \Spec A$$
for some  integer $r>0$ and $a\in A$. Here $x=(x')^\ell,\ y=(y')^\ell$. The group $\mmu_r$ acts
by the action 
$$\xi\colon(x',y')\mapsto (\xi x',\xi^m y')$$ where $m$ is an
element of $\Z\slash r$ and $\xi$ is a primitive $r$th root of the unit. The action is called \cor{balanced}
if $m\equiv -1 \mod r$. A curve with balanced action at every node is called a balanced curve.
\end{itemize}

\begin{defin}[Twisted $G$-cover]
Given an $n$-marked twisted balanced curve 
$(\Sigma_1,\dots,\Sigma_n;\ \ssC\to C\to S)$,
 a twisted $G$-cover is a representable 
stack morphism $\phi\colon\ssC\to BG$,
\cor{i.e.}~an object of the category $\Fun(\ssC,BG)$
which moreover is representable.
\end{defin}

We observe that by adding the balancing hypothesis,
we are using a slightly different notion of
twisted $G$-cover with respect to \cite{acv03}.

\begin{defin}\label{twcat}
We consider the category 
 $\mc B^{\bal}_{g,n}(G)$. 
The objects of $\mc B^{\bal}_{g,n}(G)$ are twisted $n$-marked $S$-curves of genus $g$ with a twisted
$G$-cover, for any scheme $S$. 

Consider two twisted $G$-covers $\phi'\colon \ssC'\to BG$
and $\phi\colon \ssC\to BG$ over the twisted $n$-marked curves $\ssC'$ and $\ssC$ respectively.
A morphism $(\ssC',\phi')\to(\ssC,\phi)$
is a pair $(\msf f, \alpha)$
such that $\msf f\colon \ssC'\to \ssC$
is a morphism of $n$-marked twisted curves,
and $\alpha\colon \phi'\to \phi\circ\msf f$
is an isomorphism in $\Fun(\ssC',BG)$.
\end{defin}

\subsection{Introduction to admissible $G$-covers}\label{int_adm}

\begin{defin}[Admissible cover]
Given a nodal $S$-curve $X\to S$ with marked points,
an admissible cover $u\colon F\to X$ is a morphism such that:
\begin{enumerate}
\item the composition $F\to S$ is a nodal $S$-curve;
\item given a geometric point $\bar s \in S$, every node of $F_{\bar s}$ maps via $u$ to a node of $X_{\bar s}$;
\item the restriction $F|_{X_{\gen}}\to X_{\gen}$ is an étale cover of degree $d$;
\item given a geometric point $\bar s \in S$, the local picture of $F_{\bar s}\xrightarrow{u} X_{\bar s}$ at a point of $F_{\bar s}$ mapping to a marked point of $X$
is isomorphic to
$$\Spec A[x']\to \Spec A[x]\to \Spec A,$$
for some normal ring $A$, an integer $r>0$ and $u^*x=(x')^r$;
\item the local picture of $F_{\bar s}\xrightarrow{u} X_{\bar s}$ at a node of $F_{\bar s}$
is isomorphic to
$$\Spec \left(\frac{A[x',y']}{(x'y'-a)}\right)\to \Spec \left(\frac{A[x,y]}{(xy-a^r)}\right)\to \Spec A,$$
for some integer $r>0$ and an element $a\in A$, $u^*x=(x')^r$ and $u^*y=(y')^r$.
\end{enumerate}
\end{defin}

The category $\adm_{g,n,d}$ of $n$-pointed stable curves of genus $g$
with an admissible cover of degree $d$, is a proper Deligne-Mumford stack.\newline

Consider $F\to C$ an admissible cover of a nodal curve $C$, a $G$-action on
$F$ such that the restriction $F|_{C_{\gen}}\to C_{\gen}$ is a principal $G$-bundle,
a smooth point $p$ of $C$ and a preimage $\tilde p\in F$ of $p$. We denote by $H_{\tilde p}\subset G$
the stabilizer of $\tilde p$. The $G$-action induces a primitive character
$\chi_{\tilde p}\colon H_{\tilde p}\to \GL(T_{\tilde p}F)=\C^*$.
In the set of pairs $(H,\chi)$, with $H$ a $G$ subgroup and
$\chi\colon H\to \C^*$ a character, we introduce
the equivalence relation 
$(H,\chi)\sim (H',\chi')$ if and only if there exists $s\in G$ such that $H'=sHs^{-1}$ and $\chi'=\chi^s$,
where $\chi^s\colon h\mapsto \chi(s^{-1}hs)$ is the conjugated character to $\chi$.
Consider a point $\tilde p$ on $F$ with stabilizer $H_{\tilde p}$
and associated character $\chi_{\tilde p}$. We observe that for any point $s\cdot \tilde p$ of the same fiber,
$H_{s\cdot \tilde p}=s H_{\tilde p}s^{-1}$ and $\chi_{s\cdot \tilde p}=\chi_{\tilde p}^s$.
Therefore the equivalence class of the pair $(H_{\tilde p},\chi_{\tilde p})$
only depends on the point $p$.

\begin{defin}\label{def_locgtype}
For any smooth point $\tilde p$ on $F$, we call 
\cor{local index} the associated pair $(H_{\tilde p},\chi_{\tilde p})$.
For any smooth point $p\in C$,
the conjugacy class of the local index of any $\tilde p$ in $F_p$ 
is called the \cor{$G$-type} at~$p$, following the notation in \cite{berroma11}.
We denote the $G$-type by
$\lq H_p,\chi_p\rq $,
where $H_p$ is the stabilizer of one of the points in $F_p$, and $\chi_p$ the associated character.
\end{defin}

\begin{lemma}[see {\cite[Lemma 2.3.3]{gale19}}]\label{lem_stab}
Consider $u\colon F\to C$ an admissible cover of a nodal curve $C$
such that the restriction $F|_{C_{\gen}}\to C_{\gen}$ is a principal $G$-bundle.
If $\tilde p\in F$ is one of the preimages of a node or a marked point, then the stabilizer $H_{\tilde p}$
is a cyclic group.
\end{lemma}

Observe that the set of characters $\chi\colon \mmu_r\to \C^*$ of a cyclic group, 
is the group $\Z\slash r\Z$.
In particular, the character associated to $k\in \Z\slash r\Z$ maps $\xi\mapsto \xi^k$ for any $\xi$
$r$th root of the unit. In the case of a node $\tilde p\in F$, we observe that $H_{\tilde p}$ acts independently on
the two branches $U_1$ and $U_2$. We denote by $\chi_{\tilde p}^{(1)}$ and $\chi_{\tilde p}^{(2)}$
the characters of these actions.

\begin{defin}
The $G$-action at node $\tilde p$ is balanced when $\chi_{\tilde p}^{(1)}=-\chi_{\tilde p}^{(2)}$,
that is they are opposite as elements of $\Z\slash r\Z$.
\end{defin}

\begin{defin}[Admissible $G$-cover]\label{def_admcov}
Take $G$ finite group, an admissible cover $u\colon F\to C$ of a nodal curve $C$ is an admissible $G$-cover if
\begin{enumerate}
\item the restriction $u|_{C_{\gen}}\colon F|_{C_{\gen}}\to C_{\gen}$ is a principal $G$-bundle.
This implies,
by Lemma \ref{lem_stab}, that for every node or marked point $\tilde p\in F$, the stabilizer $H_{\tilde p}$ is a cyclic group;
\item the action of $H_{\tilde p}$ is balanced for every node $\tilde p\in F$.
\end{enumerate}
\end{defin}

This notion was developed
by Abramovich, Corti and Vistoli in~\cite{acv03}, and also
by Jarvis, Kaufmann and Kimura in~\cite{jkk05}.

\begin{defin}
We call $\adm_{g,n}^G$ the stack of stable curves of genus $g$ with $n$ marked points
and equipped with an admissible $G$-cover.
\end{defin}

\begin{rmk}\label{rmk_locind}
For any cyclic subgroup $H\subset G$, we choose the privileged
root $\exp(2\pi i\slash |H|)$. After this choice, 
The datum of $(H,\chi)$, is equivalent to the datum of the $H$ generator
$h=\chi^{-1}(e^{2\pi i\slash |H|})$. As a consequence, the conjugacy class $\lq H,\chi\rq $
is identified with the conjugacy class $\lq h\rq $ of $h$ in~$G$.
\end{rmk}

\begin{defin}\label{def_hd}
Given an admissible $G$-cover $F\to C$ over an $n$-marked stable curve, the series 
$\lq h_1\rq ,\ \lq h_2\rq ,\dots,\ \lq h_n\rq$,
of the $G$-types of the singular fibers over the marked points,
is called \cor{Hurwitz datum} of the cover.
The stack of admissible $G$-covers of genus $g$ with a given
Hurwitz datum is denoted by $\adm_{g,\lq h_1\rq ,\dots,\lq h_n\rq }^G$.
\end{defin}

\begin{rmk}\label{def_gtn}
Given an admissible $G$-cover $F\to C$, if $p$ is a node of $C$
and $\tilde p$ one of its preimages on $F$, then 
the local index of $\tilde p$ and the $G$-type of $p$
are well defined once we fix a privileged branch of $p$.
Switching the branches sends the local index and the $G$-type
in their inverses.\newline
\end{rmk}

We denote by $\mc T(F)$ the set of connected components of an admissible $G$-cover $F\to C$,
and naturally it inherits the $G$-cover action $\psi$.
The conjugacy class of the stabilizer $H_E\subset G$ of any connected component $E$, is independent 
of the choice of $E$.

We call $\mc T(G)$ the set of equivalence classes in $\sub(G)$ with respect
to conjugation. Then,
for every admissible $G$-cover there exists a canonical class
$\mc H$ in $\mc T(G)$ and a canonical surjective map $\mc T(F)\twoheadrightarrow \mc H$
sending any component in its stabilizer. About the theory of group actions
that we use in this work, we also refer to \cite[\S2.1.2]{gale19}.

\begin{defin}\label{defsubcl}
Consider two subgroup conjugacy classes $\mc H_1,\mc H_2$ in $ \mc T(G)$,
we say that $\mc H_2$ is a subclass of $\mc H_1$, denoted by $\mc H_2\leq \mc H_1$,
if for one element  $H_2\in \mc H_2$ (and hence for all), there exists $H_1\in\mc H_1$
such that $H_2$ is a subgroup of $H_1$. If the inclusion is strict, then $\mc H_2$ is a
strict subclass of $\mc H_1$ and the notation is $\mc H_2<\mc H_1$.
\end{defin}

\begin{defin}\label{def_hog}
We denote by $\Ho^G(\mc T(F),G)$
the set of maps $v\colon \mc T(F)\to G$
such that $v(\psi(g,E))=g\cdot v(E)\cdot g^{-1}$.
\end{defin}

\begin{defin}\label{def_subclas}
Consider a subgroup conjugacy class $\mc H$ of $G$.
An admissible $\mc H$-cover
is an admissible $G$-cover such that every connected component has stabilizer in $\mc H$.
\end{defin}

\begin{defin}\label{def_admgh}
We denote by $\adm^{G,\mc H}_g$ the 
stack of admissible $\mc H$-covers over stable curves
of genus $g$, and we denote by
$\adm^{G,\mc H}_{g,\lq h_1\rq,\dots,\lq h_n\rq}$ the stack
of admissible $\mc H$-cover with Hurwitz datum $\lq h_1\rq,\dots,\lq h_n\rq$
over the $n$ marked points.
\end{defin}

\begin{prop}[see in {\cite[Proposition 2.3.14]{gale19}}]\label{prop_zz}
Consider $(C;p_1,\dots,p_n)$ a nodal $n$-marked curve, and
$F\to C$ an admissible $G$-cover, then
$$\Aut_{\adm}(C,F)=\Ho^G(\mc T(F), G).$$\newline
\end{prop} 

We introduced the two categories $\mc B_g^{\bal}(G)$ (see Definition \ref{twcat}) and $\adm_g^G$ with
the purpose of ``well'' defining the notion of principal $G$-bundle over stable
non-smooth curves. These two categories are proven isomorphic in \cite{acv03}.

\begin{teo}[see {\cite[Theorem 4.3.2]{acv03}}]\label{teo_3e}
There exists a base preserving equivalence between $\mc B_g^{\bal}(G)$ and $\adm_g^G$,
therefore in particular they are isomorphic Deligne-Mumford stacks.
\end{teo}

\begin{rmk}\label{rmk_not}
From now on we will use the notation $\overline\rr_{g,G}=\mc B_g^{\bal}(G)=\adm^G_g$ for the moduli stack of curves of genus $g$
equipped with an admissible $G$-cover, and $\overline\R_{g,G}$ for its coarse space.
Analogously, we will use the notation $\overline\rr_{g,G}^{\mc H}=\adm^{G,\mc H}_g$, and $\overline \R_{g,G}^{\mc H}$, for the moduli
stack of admissible $\mc H$-covers over stable curves of genus $g$, and its coarse space.
\end{rmk}

\begin{rmk}\label{rmk_equiv}
In the following, we will say that a twisted $G$-cover $(\ssC,\phi)$
``is'' an admissible $G$-cover $F\to C$ (or the other way around), meaning that $F\to C$ is
the naturally associated admissible $G$-cover to $(\ssC,\phi)$.
\end{rmk}

\begin{rmk}\label{rmk_local}
We recall the local description of $\overline\R_{g,G}$, following \cite[Remark 2.2.8]{gale19}.
For any twisted $G$-cover~$(\ssC,\phi)$, the local picture of $\overline\R_{g,G}$ at $[\ssC,\phi]$ is
$\defo(\ssC,\phi)\slash\Aut(\ssC,\phi)$,
where the action of the automorphism group is induced by the universal property. If we consider
the universal deformation $\defo(\ssC;\sing\ssC)$ of $\ssC$ alongside with its nodes, this is
naturally identified with $\defo(C;\sing C)$ where $C$ is the coarse space of $\ssC$. If $C_1,\dots, C_V$ are
the irreducible components of $C$, $\overline C_i$ their normalizations, $D_i\subset \overline C_i$ the divisors of the
the preimages of the nodes, then 
$$\defo(\ssC;\sing\ssC)=\defo(C;\sing C)=\bigoplus_{i=1}^V \defo(\overline C_i;D_i)=\bigoplus_{i=1}^VH^1(\overline C_i,T_{\overline C_i}(-D_i)).$$
We denote by $q_1,\dots,q_{\delta}$ the nodes of $C$, and we observe that we have a canonical splitting
$$\defo(\ssC)\slash\defo(\ssC;\sing\ssC)=\bigoplus_{j=1}^{\delta} R_j.$$
For every $j$, $R_j\cong\A^1$. These are the universal deformations (or smoothings) of nodes $q_j$ with the associated
stabilizers at $\ssC$. We recall that for the coarse curve $C$, we have the splitting 
$\defo(C)\slash\defo(C;\sing C)=\bigoplus_{j=1}^{\delta} M_j$,
 where
again $M_j\cong \A^1$ and there exists a canonical morphisms $R_j\to M_j$ of degree $r_j$, the cardinality of
the $q_j$ stabilizer, and branched at the origin. In particular we call $t_j$ the coordinate of $M_j$,
and $\tilde t_j$ the coordinate of $R_j$ such that $t_j=(\tilde t_j)^{r_j}$.\newline
\end{rmk}

\subsection{Monodromy description of admissible $G$-covers}\label{sec_monodromy}

Consider a smooth curve $C$ of genus $g$ and $n$ marked points $p_1,\dots, p_n$,
the fundamental group of $C_{\gen}=C\backslash\{p_1,\dots,p_n\}$
has $2g+n$ generators 
$\alpha_1,\alpha_2,\dots\alpha_g,\beta_1\dots,\beta_{g},\gamma_1,\dots,\gamma_n$.
These are represented in the figure below, where the arrows with the same label are identified
respecting the orientation. As we can also see in the figure, these generators respect the following relation,
\begin{equation}\label{eq_boh}
\alpha_1\beta_1\alpha_1^{-1}\beta_1^{-1}\cdots\alpha_{g}\beta_{g}\alpha_{g}^{-1}\beta_{g}^{-1}\cdot\gamma_1\cdots\gamma_n=1.
\end{equation}
This is called the canonical representation
of the fundamental group of a genus $g$ smooth curve.

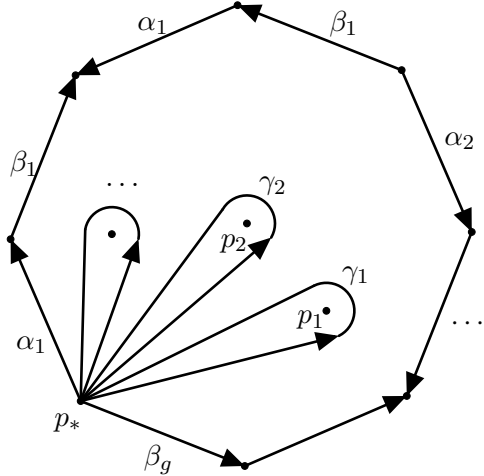
\begin{center}
\begin{tikzpicture}[scale=0.85,line cap=round,line join=round,>=triangle 45,x=1cm,y=1cm]
\draw [->,line width=1pt] (-1.93,-4.47) -- (-3.0262445840513914,-1.936568542494924);
\draw [->,line width=1pt] (-3.0262445840513914,-1.936568542494924) -- (-2.01,0.63);
\draw [->,line width=1pt] (0.5234314575050768,1.726244584051392) -- (-2.01,0.63);
\draw [->,line width=1pt] (3.09,0.71) -- (0.5234314575050769,1.726244584051392);
\draw [->,line width=1pt] (3.09,0.71) -- (4.186244584051393,-1.823431457505076);
\draw [->,line width=1pt] (4.186244584051393,-1.8234314575050763) -- (3.17,-4.39);
\draw [->,line width=1pt] (0.6365685424949233,-5.486244584051392) -- (3.17,-4.39);
\draw [->,line width=1pt] (-1.93,-4.47) -- (0.6365685424949235,-5.486244584051392);
\draw (-3.1,-3.3) node[anchor=north west] {$\alpha_1$};
\draw (-3.2,-0.4) node[anchor=north west] {$\beta_1$};
\draw (-1.2,1.7) node[anchor=north west] {$\alpha_1$};
\draw (1.8,1.8) node[anchor=north west] {$\beta_1$};
\draw (3.6,-0.1) node[anchor=north west] {$\alpha_2$};
\draw (3.7,-3.0) node[anchor=north west] {$\cdots$};
\draw (-1.1,-5) node[anchor=north west] {$\beta_g$};
\draw [shift={(-1.4427152840978867,-1.8570308865268197)},line width=1pt]  plot[domain=-0.1746721990082385:2.966920454581553,variable=\t]({1*0.4227776810856519*cos(\t r)+0*0.4227776810856519*sin(\t r)},{0*0.4227776810856519*cos(\t r)+1*0.4227776810856519*sin(\t r)});
\draw [shift={(0.67,-1.69)},line width=1pt]  plot[domain=-0.5247957716501066:2.616796881939686,variable=\t]({1*0.4390899680020028*cos(\t r)+0*0.4390899680020028*sin(\t r)},{0*0.4390899680020028*cos(\t r)+1*0.4390899680020028*sin(\t r)});
\draw [shift={(1.912531644340949,-3.057082707647277)},line width=1pt]  plot[domain=-1.1071487177940948:2.0344439357957063,variable=\t]({1*0.43810570588534736*cos(\t r)+0*0.43810570588534736*sin(\t r)},{0*0.43810570588534736*cos(\t r)+1*0.43810570588534736*sin(\t r)});
\draw [->,line width=1pt] (-1.93,-4.47) -- (-1.0263707747295636,-1.930503447003582);
\draw [->,line width=1pt] (-1.93,-4.47) -- (1.05,-1.91);
\draw [->,line width=1pt] (-1.93,-4.47) -- (2.108458472278981,-3.4489363635233445);
\draw [line width=1pt] (-1.8590597934662094,-1.7835583260500565)-- (-1.93,-4.47);
\draw [line width=1pt] (0.29,-1.47)-- (-1.93,-4.47);
\draw [line width=1pt] (1.7166048164029144,-2.665229051771211)-- (-1.93,-4.47);
\draw (-1.7,-0.9) node[anchor=north west] {$\dots$};
\draw (0.7,-0.8) node[anchor=north west] {$\gamma_2$};
\draw (2,-2.2) node[anchor=north west] {$\gamma_1$};
\draw (0.1,-1.7) node[anchor=north west] {$p_2$};
\draw (1.3,-2.9) node[anchor=north west] {$p_1$};
\draw (-2.5,-4.5) node[anchor=north west] {$p_*$};
\begin{scriptsize}
\draw [fill=black] (-1.93,-4.47) circle (1.5pt);
\draw [fill=black] (3.17,-4.39) circle (1.5pt);
\draw [fill=black] (3.09,0.71) circle (1.5pt);
\draw [fill=black] (-2.01,0.63) circle (1.5pt);
\draw [fill=black] (-3.0262445840513914,-1.936568542494924) circle (1.5pt);
\draw [fill=black] (0.5234314575050768,1.726244584051392) circle (1.5pt);
\draw [fill=black] (4.186244584051393,-1.8234314575050763) circle (1.5pt);
\draw [fill=black] (0.6365685424949233,-5.486244584051392) circle (1.5pt);
\draw [fill=black] (-1.4427152840978867,-1.8570308865268192) circle (1.5pt);
\draw [fill=black] (0.67,-1.69) circle (1.5pt);
\draw [fill=black] (1.9125316443409486,-3.0570827076472797) circle (1.5pt);
\end{scriptsize}
\end{tikzpicture}
\end{center}

It is possible to describe admissible $G$-covers over smooth curves
by the monodromy action, as done for example in \cite[\S 2.3]{berroma11} and \cite[\S 3.5]{schmizel18}.
Consider a smooth curve $C$, a generic point $p_*$ on it
and an admissible $G$-cover $F\to C$. We denote the points of the fiber $F_{p_*}$ by
$\tilde p_*^{(g)}$ for any $g\in G$, in such a way that $g\cdot \tilde p_*^{(1)}=\tilde p_*^{(g)}$.
This induces a group morphism $\pi_1(C_{\gen},p_*)\to G$.
This monodromy morphism is well defined
up to relabelling  the points $\tilde p_*^{(g)}$, \cor{i.e.}~up to
$G$ conjugation.
The following proposition is a rephrasing of \cite[Lemma 2.6]{berroma11}.

\begin{prop}\label{prop_gfg2}
Given a smooth $n$-marked curve $(C;p_1,\dots,p_n)$
and a point $p_*$ on its generic locus $C_{\gen}$,
the set of isomorphism classes of admissible $G$-covers on $C$
 is naturally in bijection
with the set of conjugacy classes of maps 
$$\varpi\colon \pi_1(C_{\gen},p_*)\to G.$$
\end{prop}
\begin{rmk}\label{rmk_mono}
We also point out that, as represented in the figure, the monodromy of $\gamma_i$
 at any point~$p_*^{(g)}$, with $g\in G$, is given by a small circular
lacet around the deleted point~$p_i$. Therefore
by definition of $G$-type, if $\lq h_i\rq $ is the $G$-type of $p_i$,
then $\lq \varpi(\gamma_i)\rq=\lq h_i\rq$.
\end{rmk}

We consider now the case of an irreducible curve $C$ with one type $0$ node, \cor{i.e.}~an autointersection node.
If $\overline C$ is the $C$ normalization, we denote by $p_1$ and $p_2$ the preimages of the node on $\overline C$.
An admissible $G$-cover $F\to C$ induces an admissible $G$-cover $\overline F\to \overline C$ on the $2$-marked
genus $g-1$ curve $(\overline C;p_1,p_2)$ and by what we said above, there exists an associated
morphism $\varpi\colon\pi_1(\overline C_{\gen})\to G$ such that $\lq \varpi(\gamma_1)\rq=\lq\varpi(\gamma_2)^{-1}\rq$.
As we obtain $C$ from $\overline C$ by gluing $p_1$ to $p_2$, $F$ is obtained by gluing the fibers $\overline F_{p_1}$ and $\overline F_{p_2}$.

\begin{rmk}\label{mono2}
In the last case, consider a path $\gamma_{1,2}$ on $\overline C$ from $p_1$ to $p_2$. Gluing $\overline F_{p_1}$ to $\overline F_{p_2}$
is equivalent to lift $\gamma_{1,2}$ and therefore to give a monodromy factor $h_{\gamma}$ well defined up to conjugation.
Taking the twisted $G$-cover point of view, a point of $\overline F_{p_1}$ has to be sent by $h_{\gamma}$ to a point
with the same local index (see \cite[\S2.2.2]{gale19}).  For example if we take a point with local index
$w\in\lq\varpi(\gamma_1)\rq$,
then the local index of its image is $h_\gamma w h_\gamma^{-1}$.
This means that $w$ and
$h_{\gamma}$ must commute. This condition allows to define $h_\gamma$ up to conjugation and is also sufficient to make the gluing.
\end{rmk}

\begin{rmk}\label{rmk_autoin}
If $C$ is a nodal curve with a type $0$ node as above, the fundamental group of the normalization
$\pi_1(\overline C_{\gen})$ is naturally a subgroup of $\pi_1(C)$, therefore the generators $\alpha_1,\dots,\beta_{g-1}$
are $\pi_1(C)$ elements.
Consider a slight deformation $F'\to C'$ of $F\to C$, smoothing the node, such that
$C'$ is a smooth curve. This induces a natural morphism $\pi_1(C)\to \pi_1(C')$, in particular the $\alpha_i,\beta_i$
 from $\pi_1(\overline C_{\gen})$ are generators in the canonical representation of $\pi_1(C')$,
and we can complete this representation in such a way that the class of lacet $\gamma_1$ becomes the element $\alpha_g$ in $\pi_1(C')$.
The generator $\beta_g$ comes from the gluing data of the nodal fiber.
\end{rmk}

Following this monodromy approach, and
in order to describe the components of $\overline\rr_{g,S_3}$ in~\S\ref{sec_RRR},
it is useful to introduce a notion of group morphism defined up to conjugation.

\begin{defin}\label{def_surj}
Consider two groups $A,\ G$, and the equivalence relation $\sim$
such that for two morphisms $\varphi,\varphi'\colon A\to G$, $\varphi\sim \varphi'$
if there exists $g\in G$, $\varphi'=g\cdot \varphi\cdot g^{-1}$.

For any conjugacy class $\mc H\in \mc T(G)$,
we denote by $\Ho_{\m{s}}(A,\mc H)$ the subset of $\Ho(A,G)\slash \sim$
induced by the morphisms whose image is a subgroup in $\mc H$.

We will use the notation $\varphi\colon A\overset{G}{\rightsquigarrow} \mc H$ for 
an element of $\Ho_{\m{s}}(A,\mc H)$, that is a group morphism $A\to G$
defined up to conjugation whose image is in $\mc H$.\newline
\end{defin}

\section{Structure of $\overline{\mathcal{R}}_{g,S_3}$}\label{struct}

In the first section we study the boundary of $\overline\R_{g,S_3}$,
\cor{i.e.}~the locus $\overline\R_{g,G}\backslash\R_{g,G}$. 
In the other two sections
we describe the canonical divisor of this space as a combination
of the Hodge divisor
and boundary divisors.

\subsection{Moduli boundary of $\overline{\mathcal{R}}_{g,S_3}$}

For every subgroup conjugacy class $\mc H\in\mc T(G)$, we 
consider the moduli 
stack $\overline\rr_{g,G}^{\mc H}$
of admissible $\mc H$-covers over stable curves of genus $g$, \cor{i.e.}~admissible
$G$-covers such that the stabilizer of every connected component of a cover is in 
the class $\mc H$ (see Definition \ref{def_subclas} and Remark \ref{rmk_not}). 
These are all pairwise disconnected substacks of~$\overline\rr_{g,G}$.

\subsubsection{Components of $\overline\R_{g,S_3}$}\label{sec_RRR}

We observe that $S_3$ has $6$ subgroups 
and $4$ subgroup classes.

\begin{itemize}
\item The trivial subgroup $(1)\subset S_3$. In this case $\overline\R_{g,S_3}^{1}$
is isomorphic to $\overline\M_g$.
\item The three subgroups $T_1,T_2,T_3$ of order $2$, generated respectively
by transpositions $(23),(13),(12)$. These subgroups all stay in the same conjugacy class $T=\{T_1,T_2,T_3\}$.
We are going to show that the stack $\rr_{g,S_3}^{T}$ is isomorphic to $\rr_{g,\mmu_2}^{\mmu_2}\subset\rr_{g,2}^0$,
the moduli stack of twisted curves equipped with a non-trivial square root of the trivial bundle.

Indeed, any admissible connected $\mmu_2$-cover $E\to C$ is equivalent to the data $(C,\rho)$
of a curve $C$ and $\rho\in\Ho_{\m{s}}(\pi_1(C),\mmu_2)$, as stated by Proposition \ref{prop_gfg2} and
using the notation of Definition \ref{def_surj}.
There exists an isomorphism $\sigma_2\colon\mmu_2\overset{S_3}{\rightsquigarrow} T$ (see again Definition \ref{def_surj}),
and this induces a set bijection
$$\Ho_{\m{s}}(\pi_1(C),\mmu_2)\cong \Ho_{\m{s}}(\pi_1(C), T).$$ 
Therefore the map $\rr_{g,\mmu_2}^{\mmu_2}\to\rr_{g,S_3}^T$ defined by
$(C,\rho) \mapsto (C,\sigma_2\circ\rho)$ is an isomorphism.
We additionally observe that given an admissible $T$-cover $F\to C$ and any connected component $E\subset F$,
the map $E\to C$ is 
an admissible connected $\mmu_2$-cover.

We observe that the isomorphism of moduli stacks extends to $\overline\rr_{g,S_3}^T\cong\overline\rr_{g,\mmu_2}^{\mmu_2}$,
but the construction at the boundary introduces an additional complexity which is useless for the purpose of this work.

\item The normal subgroup $N\subset S_3$, a cyclic group generated by
the $3$-cycle $(123)$. With a little abuse of notation
we call $N$ the class $\{N\}\in\mc T(S_3) $.
Consider the moduli stack $ \rr_{g,\mmu_3}^{\mmu_3}\subset\rr_{g,\mmu_3}^0$
of twisted curves equipped with a non-trivial
third root of the trivial bundle. We are going to prove that there exists a $2:1$ map $\rr_{g,\mmu_3}^{\mmu_3}\to \rr_{g,S_3}^N$.

Indeed, any admissible connected $\mmu_3$-cover $E\to C$ is equivalent to the data $(C,\rho)$
of a curve $C$ and an element $\rho\in\Ho_{\m{s}}(\pi_1(C),\mmu_3)$. There exists an
isomorphism $\sigma_3\colon\mmu_3\overset{S_3}{\rightsquigarrow} N$,
and this induces a set surjection
$$\Ho_{\m{s}}(\pi_1(C),\mmu_3)\twoheadrightarrow \Ho_{\m{s}}(\pi_1(C),N).$$
If $\we 2\colon\mmu_3\to\mmu_3$ is the second power map, then $\we 2$ also acts as an involution on $\Ho_{\m{s}}(\pi_1(C),\mmu_3)$.
We observe that
if $\rho\colon\pi_1(C)\overset{\mmu_3}{\rightsquigarrow} \mmu_3$, then $\sigma_3\circ\rho=\sigma_3\circ(\we 2)\circ\rho$. 
In fact, the following is true
$$\Ho_{\m{s}}(\pi_1(C),\mmu_3)\slash(\we 2)\cong \Ho_{\m{s}}(\pi_1(C),N).$$
Therefore the map $(C,\rho)\mapsto(C,\sigma_3\circ\rho)$ is the $2:1$ map we were searching. As before
the map extends to $\overline\rr_{g,\mmu_3}^{\mmu_3}\to\overline\rr_{g,S_3}^N$.

\item The group $S_3$ itself. In this case the stack $\overline\R_{g,S_3}^{S_3}$ is the moduli
of curves equipped with a connected admissible  $S_3$-cover. This is the
``really new'' component of moduli space $\overline\R_{g,S_3}$,
and our analysis will
focus on it.
\end{itemize}

Furthermore, we observe that there exists a canonical map $\overline\rr_{g,S_3}\to \overline\rr_{g,\mmu_2}$.
Any admissible $G$-cover $F\to C$
is equivalent to the data of a curve $C$ plus an element of $\Ho_{\m{s}}(\pi_1(C),S_3)$.
The quotient $S_3\twoheadrightarrow S_3\slash N=\mmu_2$ induces a surjection
$\Ho_{\m{s}}(\pi_1(C),S_3)\twoheadrightarrow \Ho_{\m{s}}(\pi_1(C),\mmu_2)$ and therefore the map above.\newline

\subsubsection{The boundary divisors}
To classify the boundary divisors of $\overline\R_{g,G}$, that is the divisors filling the locus $\overline\R_{g,G}\backslash \R_{g,G}$, we start
by recalling the boundary divisors of $\overline\M_g$.

\begin{defin}
A node $q$ of a nodal curve $C$ is a disconnecting node if the partial
normalization of $C$ at $q$ is a disconnected scheme. This disconnected scheme
has two connected components noted $C_1$ and $C_2$.
Consider $C$ of genus $g$ and $q$ a disconnecting node,
$q$ is of type~$i$, with $1\leq i\leq\lfloor g\slash 2\rfloor$, if $C_1$ and $C_2$ 
have genus $i$ and $g-i$.
If $q$ is a non-disconnecting node, then it is called a node of type $0$.
\end{defin}

\begin{defin}
For every $i$ with $1\leq i\leq \lfloor g\slash 2\rfloor$, the divisor $\Delta_i\subset\overline\M_g$
is the locus of curves with a disconnecting node of type $i$.
For $i=0$, the divisor $\Delta_0\subset\overline\M_g$ is the locus of curves
with a node of type $0$, or equivalently the closure of the locus of
nodal irreducible curves.
\end{defin}

For every $i$ we denote by $\delta_i$ the class of $\Delta_i$ in $\Pic_\Q(\overline\M_g)$.
We consider the natural morphism $\pi\colon\overline\R_{g,G}\to \overline\M_g$ and 
look at the preimages $\pi^{-1}(\Delta_i)$ for every $i$. The intersections of this preimages
with the connected components of $\overline\R_{g,G}$ are the boundary components that we are going to consider.

We start by focusing on the loci of curves with a disconnecting node, \cor{i.e.}~the preimages
of $\Delta_i$ with $i\neq 0$. In what follows for any curve $C$ in $\Delta_i$ we denote by
$q$ the node who separates the components $C_1$, of genus $i$, and $C_2$, of genus $g-i$.
As a consequence
$C_1\sqcup C_2 \to C$
is the partial normalization of $C$ at $q$, and we denote by $q_1$ and $q_2$ the preimages of $q$
respectively on $C_1$ and $C_2$.

There exists, at the level of the moduli space $\overline\M_g$, a natural gluing map
$$\overline\M_{i,1}\times\overline\M_{g-i,1}\to\overline\M_g,$$
defined in such a way that a pair of points $([C_1,q_1],[C_2,q_2])$, where
$C_1$ and $C_2$ have genus $i$ and $g-i$ respectively,
is sent on the point associated to the nodal curve 
$$C:=(C_1\sqcup C_2)\slash(q_1\sim q_2).$$
We recall that any admissible $G$-cover on $C$
induces two admissible $G$-covers $F_1$ on $(C_1;q_1)$ and $F_2$ on $(C_2;q_2)$
such that the $G$-types $\lq h_1\rq $ and $\lq h_2\rq $ on $q_1$ and $q_2$
are one the inverse of the other, $\lq h_1\rq =\lq h_2^{-1}\rq $.

\begin{defin}
Consider two subgroup conjugacy classes $\mc H_1 $ and $\mc H_2$ in~$\mc T(G)$,	
a conjugacy class $\lq h\rq $ in $\lq G\rq$ and $i$ such that  $1\leq i \leq \lfloor g\slash 2\rfloor$.
We denote by
$\Delta_{i,\lq h\rq }^{\mc H_1,\mc H_2 }$
the locus in $\overline\R_{g,G}$ of curves $C$ with a node $q$ of type $i$ and
with an admissible $G$-cover $F\to C$ such that if we call 
$F_1\to (C_1;q_1)$ and $F_2\to (C_2;q_2)$ the restrictions of~$F$,
then these are an admissible $\mc H_1$-cover of $(C_1;q_1)$ and an admissible $\mc H_2$-cover
of $(C_2;q_2)$
respectively (see Definition \ref{def_subclas}), and the $G$-type at $q$, with respect to the branch of $C_1$, is~$\lq h\rq$.
In order to simplify the notation we will omit the $G$-type when it is trivial.
\end{defin}

The loci $\Delta_{i,\lq h\rq }^{\mc H_1,\mc H_2 }$ are (not necessarily connected) 
divisors of the moduli space~$\overline\R_{g,G}$.

For example
in the case of $\overline\R_{g,S_3}$ we have 
the following classes:
\begin{itemize}
\item $\Delta_i^{1,S_3},\ \Delta_i^{S_3,1},\ \Delta_i^{1,T},\ \Delta_i^{T,1},\ \Delta_i^{1,N},\ \Delta_i^{N,1}$ and $\Delta_i^{1,1}\cong\Delta_i$
are the cases of admissible $G$-covers
which are trivial over $C_1$ or $C_2$;
\item $\Delta_i^{T,T},\Delta_i^{T,N},\Delta_i^{N,T},\Delta_i^{T,S_3},\Delta_i^{N,N},\Delta_i^{N,S_3},\Delta_i^{S_3,N}$ and $\Delta_i^{S_3,S_3}$
are the other cases with trivial $G$-type at the node $q$;
\item $\Delta_{i,c_3}^{S_3,S_3}$ is the only case with non-trivial stabilizer at $q$. Here we denoted by $c_3$ the conjugacy class
$\lq 123\rq=\{(123),(132)\}$.
\end{itemize}

\begin{rmk}
We observe that by Equation (\ref{eq_boh}) and Remark \ref{rmk_mono}, it is possible to have non-trivial stabilizer
at $q$ only if for both covers $F_1\to C_1$ and $F_2\to C_2$,
the associated morphisms $\pi_1(C_i)\to S_3$ are surjective, and therefore
both are connected admissible $S_3$-covers.\newline
\end{rmk}

Finally, we consider the preimage of $\Delta_0$, \cor{i.e.}~the locus of curves with a node of type~$0$.
We start by working on a covering stack of $\pi^{-1}(\Delta_0)$.
\begin{defin}
Category $\msf D_0$ has for objects the data of a curve $C$ with a node $q$ of type $0$
with an admissible $G$-cover $F\to C$ and a privileged branch at $q$. Morphisms are admissible $G$-cover
morphisms preserving the privileged branch.
\end{defin}
\begin{rmk}
Category $\msf D_0$ is a Deligne-Mumford stack. Its coarse space $\D_0$ comes with a natural $2:1$
morphism $\gamma\colon\D_0\to\pi^{-1}(\Delta_0)$.
\end{rmk}
In what follows we consider a curve $C$
with a node $q$ of type $0$, we call 
$\nor\colon\overline C\to C$
the partial normalization of $C$ at $q$.
Given any admissible $G$-cover $F\to C$, the pullback
$\nor^*F$ over $\overline C$ is still an admissible $G$-cover.

\begin{defin}
Consider two conjugacy classes $\mc H_1,\mc H_2$ in $\mc T(G)$ such that
$\mc H_2\leq \mc H_1$ (see Definition \ref{defsubcl}).
The category $\msf D_{0,\lq h\rq }^{\mc H_1,\mc H_2}$ is a full subcategory of $\msf D_0$.
Its objects are stable curves $C$ with a node $q$ of type $0$, a privileged branch at $q$
and an admissible $\mc H_1$-cover $F\to C$
such that the pullback $\nor^*F\to\overline C$ is an admissible $\mc H_2$-cover on $\overline C$, and the 
$G$-type of $F$ at $q$ with respect to the privileged branch is $\lq h\rq $.
In case $h=1$ we will omit its notation.
The category $\msf D_{0,\lq h\rq }^{\mc H_1,\mc H_2}$ is again a Deligne-Mumford stack and we denote
its coarse space by~$\D_{0,\lq h\rq }^{\mc H_1,\mc H_2}$.
\end{defin}

\begin{rmk}
By Equation (\ref{eq_boh}) and Remark \ref{rmk_mono},
there exists a compatibility condition for the stack $\msf D_{0,\lq h\rq }^{\mc H_1,\mc H_2}$
to be non-empty. There must exist an element $h$ in the class $\lq h\rq $ and a subgroup $H_2$ in
the class $\mc H_2$ such that $h$ lies in the subgroup of commutators of~$H_2$.
\end{rmk}

\begin{rmk}
There exists a natural automorphism
$\inv\colon\msf D_0\to \msf D_0$,
which sends any curve with an admissibile $G$-cover
to the same curve and $G$-cover but changing the privileged branch at the node.
This sends isomorphically $\msf D_{0,\lq h\rq }^{\mc H_1,\mc H_2 }$ in~$\msf D_{0,\lq h^{-1}\rq }^{\mc H_1,\mc H_2}$.
\end{rmk}

\begin{defin}
For any group $G$ we consider the inverse relation
in the set of conjugacy classes $\lq G\rq $, that is
$\lq h\rq \sim \lq h'\rq $
if and only if $\lq h\rq =\lq h'\rq $ or $\lq h\rq =\lq h'^{-1}\rq $. Then we define the set $\underline{ \lq G\rq }$ as
$\underline{\lq G\rq }:=\lq G\rq \slash \sim$.
we denote by $\underline{\lq h\rq} $ the class in $\underline{\lq G\rq }$ of any
element~$h\in G$.
\end{defin}

We observe that any point $[C,F]$ of $\D_0$ and its image $\inv([C,F])$
are sent by $\gamma\colon\D_0\to\pi^{-1}(\Delta_0)$ to the same point.
This means that the image of $\D_{0,\lq h\rq }^{\mc H_1,\mc H_2 }$ via $\gamma$
depends only on $\mc H_1,\mc H_2 $ and on the class $\underline{\lq h\rq} \in\underline{\lq G\rq }$.
Equivalently, for all ${\lq h\rq} \in{\lq G\rq }$, $\D_{0,\lq h\rq }^{\mc H_1,\mc H_2 }$ and $\D_{0,\lq h^{-1}\rq }^{\mc H_1,\mc H_2}$
have the same image.

\begin{defin}
For every class $\mc H_1,\mc H_2 $ in $\mc T(G)$ and $\underline{\lq h\rq} $ in $\underline {\lq G\rq }$,
the locus $\Delta_{0,\underline{\lq   h\rq} }^{\mc H_1,\mc H_2 }$ is the image of $\D_{0,\underline{\lq h\rq} }^{\mc H_1,\mc H_2 }$ via $\gamma$
for any $G$-type $\lq h\rq\in\underline{\lq h\rq}$. As before we will omit to note $\underline{\lq h\rq} $ in the case
of the trivial class $\underline{\lq 1\rq} $.
\end{defin}

As before, we list the divisors in $\pi^{-1}(\Delta_0)$:
\begin{itemize}
\item $\Delta_0^{N,1},\ \Delta_0^{T,1}$ and $\Delta_0^{1,1}\cong\Delta_0$
are the divisor of admissible $G$-covers $F\to C$ such that $\nor^*F\to \overline C$ is trivial;
\item $\Delta_0^{T,T},\ \Delta_0^{S_3,T},\ \Delta_0^{N,N}, \Delta_0^{S_3,N}$ and $\Delta_0^{S_3,S_3}$
are the other cases with trivial associated $G$-type at node $q$;
\item with another small abuse of notation, we define $c_2:=\underline {\lq 12\rq}=\{(12),(13),(23)\}$
and $c_3:=\underline{\lq 123\rq}=\{(123),(132)\}$. Then, 
$\Delta_{0,c_2}^{T,T},\ \Delta_{0,c_2}^{S_3,S_3},\ \Delta_{0,c_3}^{N,N}$ and $\Delta_{0,c_3}^{S_3,S_3}$
are the cases with non-trivial stabilizer at $q$.
\end{itemize}

Observe that
the divisors $\Delta_{0,c_2}^{S_3,T}$ and $\Delta_{0,c_3}^{S_3,N}$ are both empty
as a direct consequence of Remark \ref{rmk_autoin}.\newline

\subsection{The canonical divisor}

To evaluate the canonical divisor we start by evaluating the pullbacks of the $\overline\M_g$ boundary divisors.
We denote by $\delta_{i,\lq h\rq }^{\mc H_1,\mc H_2 }$ the class of the divisor $\Delta_{i,\lq h\rq }^{\mc H_1,\mc H_2 }$ in the
ring~$\Pic_\Q(\overline\R_{g,S_3})$. Furthermore, we introduce some other notation to simplify the formulas
that will follow.
\begin{itemize}
\item For any $i\geq1$ we define
$$\delta_i':=\sum \delta_i^{	\mc H_1,\mc H_2},$$
where the sum is over al the divisors with trivial $S_3$-type at the node.
\item For $i=0$ we define
\begin{align*}
\delta_0'&:=\sum_{\mc H_2\leq \mc H_1}\delta_0^{\mc H_1,\mc H_2},\\
\delta_{0,c_2}&:=\delta_{0,c_2}^{T,T}+\delta_{0,c_2}^{S_3,S_3},\\
\delta_{0,c_3}&:=\delta_{0,c_3}^{N,N}+\delta_{0,c_3}^{S_3,S_3}.
\end{align*}
\end{itemize}

Given the classes $\delta_i$ in $\Pic_\Q(\overline\M_g)$, and the natural morphism
$\pi\colon\overline\R_{g,S_3}\to \overline\M_g$, we consider the pullbacks $\pi^*\delta_i$
to the ring $\Pic_\Q(\overline\R_{g,S_3})$.
\begin{lemma}\label{lemmai2}
If $i>1$, then
$$\pi^*\delta_i=\delta_i'+3\delta_{i,c_3}^{S_3,S_3}.$$
\end{lemma}
\proof This relation is true also in $\Pic_\Q(\overline\rr_{g,S_3})$.
We consider a general $S_3$-cover $(\ssC,\phi)$ in the divisor $\Delta_{i,c_3}^{S_3,S_3}$,
and we denote by $C$ the coarse space of the curve and by $q$ the type $i$ node. By construction we have
$$\frac{\Aut(\ssC,\phi)}{\Aut(C)}=\Aut(\ssC_q)=\mmu_3.$$
This implies that the morphism $\overline\rr_{g,S_3}\to\overline\mm_g$ is $3$-ramified
over $\Delta_{i,c_3}^{S_3,S_3}$. Over $\Delta_i$ but outside $\Delta_{i,c_3}^{S_3,S_3}$ the morphism
is étale, and the thesis follows.\fine

\begin{lemma}\label{lemmai1}
If $i=1$, then
$$\pi^*\delta_1=\delta_1'+3\delta_{1,c_3}^{S_3,S_3}.$$
\end{lemma}
\proof The case of a curve $\ssC$ with an elliptic tail, is different from the case $i>1$ because
there exists an elliptic
tail involution on $\ssC$. This involution
always lifts to an admissible $S_3$-cover.
Indeed, an admissible 
$S_3$-cover over an elliptic curve $(E,p_*)$ is the datum
of a morphism $\pi_1(E,p_*)\to S_3$ defined up to conjugation, and
this is always preserved by the involution.
Thus, $\pi^*\Delta_1=\Delta_1'+\Delta_{1,c_3}^{S_3,S_3}$.
To obtain the result we observe that, as in the case of Lemma \ref{lemmai2},
the morphism $\overline\rr_{g,S_3}\to\overline\mm_g$ is $3$-ramified
over $\Delta_{1,c_3}^{S_3,S_3}$.\fine

\begin{lemma}\label{lemmai0}
If $i=0$, then
$$\pi^*\delta_0=\delta_0'+2\delta_{0,c_2}+3\delta_{0,c_3}.$$
\end{lemma}
\proof Similarly to what observed in the case of Lemma \ref{lemmai2},
the morphism $\overline\rr_{g,S_3}\to\overline\mm_g$ is $2$-ramified over $\Delta_{0,c_2}$,
 $3$-ramified over $\Delta_{0,c_3}$ and étale elsewhere over $\Delta_0$.\fine

From the Harris and Mumford work \cite{harmum82}, we know the evaluation 
of the canonical divisor $K_{\overline\M_g}$ of the moduli space $\overline\M_g$.
Knowing that the morphism $\pi\colon\overline\R_{g,S_3}\to\overline\M_g$ is étale
outside the divisors $\Delta_i$ and a sublocus of codimension greater than $1$, we can
now evaluate the canonical divisor $K_{\overline\R_{g,S_3}}$.

\begin{lemma}\label{lemma_kk}
On the smooth variety $\overline\R_{g,S_3}^{\reg}$, the sublocus of regular points,
we have the following evaluation of the canonical divisor,
$$K_{\overline\R_{g,S_3}}=13\lambda-(2\delta_0'+3\delta_{0,c_2}+4\delta_{0,c_3})
-(3\delta_1'+7\delta_{1,c_3}^{S_3,S_3})-\sum_{i=2}^{\lfloor g\slash 2\rfloor}(2\delta_i'+4\delta_{i,c_3}^{S_3,S_3}).$$
\end{lemma}

\proof As proved in \cite[Theorem 2 bis, p.52]{harmum82},
on $\overline\M_g^{\reg}$ we have
$$K_{\overline\M_g}=13\lambda-2\delta_0-3\delta_1-2\delta_2-\cdots-2\delta_{\lfloor g\slash2 \rfloor}.$$
We know, from the description above, that the ramification divisor
of $\pi\colon\overline\R_{g,S_3}\to \overline\M_g$ is
 $$R=\delta_{0,c_2}+2\delta_{0,c_3}+2\sum_{i=1}^{\lfloor g\slash 2\rfloor}\delta_{i,c_3}^{S_3,S_3}.$$
By Hurwitz formula we have 
$K_{\overline\R_{g,S_3}}=\pi^*K_{\overline\M_g}+R$,
this is true on $\pi^{-1}(\overline\M_g^{\reg})$, but this locus and $\overline\R_{g,S_3}^{\reg}$
differ by a codimension $2$ locus, therefore the evaluation
is unchanged on $\overline\R_{g,S_3}^{\reg}$. Precisely, if
$\pi$ is non-étale in a point $[C,F]\in\overline\R_{g,S_3}$ outside the $\pi^{-1}(\Delta_i)$, 
then necessarily $C$ is an irreducible curve with a non-trivial automorphism group.
The thesis follows by the evaluations of the divisors $\pi^*\delta_i$ of Lemmata
\ref{lemmai2}, \ref{lemmai1} and \ref{lemmai0}.\fine

\subsection{The subspace of covers on irreducible curves}
We consider the full substack  $\widetilde{\rr}_{g,S_3}$ of $\overline{\rr}_{g,S_3}$
of irreducible stable curves of genus $g$ with an admissible $S_3$-cover. 
We denote by $\widetilde\R_{g,S_3}$ its coarse space.
We will prove that the bigness of the canonical divisor $K_{\overline\R_{g,S_3}}$ over~$\widetilde\R_{g,S_3}$
implies the bigness over~$\overline\R_{g,S_3}$. As a consequence,
it suffices to prove that $K_{\overline\R_{g,S_3}}$ is big
over~$\widetilde\R_{g,S_3}$, to prove $\overline\R_{g,S_3}$
being of general type.

We build some pencils filling up specific divisors of $\overline\R_{g,S_3}^{S_3}$.
To do this we follow \cite{farlud10} and~\cite{farpop05}. 
Given a general K3 surface $X$ of degree $2i-2$ in $\p^i$, the map
$\Bl_{i^2}(X)\to\p^1$
is a family of genus $i$ stable curves lying on $X$,
where $\Bl_{i^2}(X)$ is the blowup of $X$ in $i^2$ points. This induces a pencil $B\subset \overline\M_i$.
Moreover, there exists at least one section $\sigma$ on $B$, therefore for
every genus $g>i$ we can glue along $\sigma$
a fixed $1$-marked curve $(C_2,p_*)$ of genus $g-i$, thus inducing a pencil
$B_i\subset\overline\M_g$. The pencils $B_i$ fill up the divisors $\Delta_i$
except for $i=10$, if $i=10$ then the $B_{10}$ fill up a divisor $\mc Z$ in $\Delta_{10}$
which is the locus of smooth curves of genus $10$ lying on a K3 surface and attached
to a curve of genus $g-10$. We denote by $\mc Z'$ the preimage 
of $\mc Z$ in $\overline\R_{g,S_3}^{S_3}$.

\begin{lemma}\label{lemma_eff}
Consider an effective divisor $E$ in $\Pic_\Q(\R_{g,S_3}^{S_3})$ such that
$$[\overline E]=a\cdot\lambda -b_0'\cdot\delta_0'-b_{0,c_2}\cdot\delta_{0,c_2}-b_{0,c_3}\cdot \delta_{0,c_3}-
b_1'\cdot\delta_1'
-b_{1,c_3}^{S_3,S_3}\cdot \delta_{1,c_3}^{S_3,S_3}
-\sum_{i=2}^{\lfloor g\slash 2\rfloor}(b_i'\cdot \delta_i'+b_{i,c_3}^{S_3,S_3}\cdot\delta_{i,c_3}^{S_3,S_3})$$
in  $\Pic_\Q(\overline\R_{g,S_3}^{S_3})$,
where the $a$ and $b$ are rational coefficients.

If $a\leq 13$, $b_0'\geq 2$, $b_{0,c_2}\geq 3$  and $b_{0,c_3}\geq 4$, then $b_i'\geq 3$ and
 $b_{i,c_3}^{S_3,S_3}> 7$
for all $i\geq 1$ and~$i\neq 10$. The same is true
for $i=10$ if $\overline E$ does not contain $\mc Z'$.
\end{lemma}

\proof  By \cite[Lemma 2.4]{farpop05} we have
$$B_i\cdot \lambda=i+1,\ \ B_i\cdot \delta_0= 6i+18,\ \ B_i\cdot \delta_i=-1,\ \ B_i\cdot \delta_j=0\ \forall j\neq 0,i.$$
We introduce some pencils lying in the preimage of $B_i$ with respect
to the natural projection $\pi\colon\overline\R_{g,S_3}^{S_3}\to\overline\M_g$.

The pencil $A_i^{T,N}$ is
the preimage of $B_i$ in $\Delta_i^{T,N}$. It is
 obtained
by taking any $1$-marked stable curve $C_1$ of genus $i$ with an admissible $T$-cover, 
and gluing it to a fixed general $1$-marked curve $C_2$ of genus $g-i$ with an admissible $N$-cover.
We remark that the gluing is uniquely defined.

The pencil $A_{i,c_3}^{S_3,S_3}$ is the preimage of $B_i$ in $\Delta_{i,c_3}^{S_3,S_3}$.
It is obtained by taking any $1$-marked stable curve $(C_1,p_*')$ of genus $i$ with an admissible
connected $S_3$-cover and $S_3$-type equal to $c_3$ at $p_*'$,
and gluing it to a fixed general $1$-marked curve $(C_2,p_*)$ of genus $g-i$
with an admissible connected $S_3$ cover and the same $S_3$-type at $p_*$.
Again, the gluing is unique.\newline

We can write down some intersection numbers for the $A_i^{T,N}$:

\begin{itemize}
\item $A_i^{T,N}\cdot \lambda=(2^{2i}-1)(i+1)$, which is true because $\pi^*\lambda_{\overline\M_g}=\lambda_{\overline\R_{g,S_3}}$
and moreover $\pi_*A_i^{T,N}=B_i\cdot \deg(A_i^{T,N}\slash B_i)$;
\item $A_i^{T,N}\cdot\delta_0'=(2^{2i-1}+1)(6i+18),\ A_i^{T,N}\cdot \delta_{0,c_2}=(2^{2i-1}-2)(6i+18),\ A_i^{T,N}\cdot\delta_{0,c_3}=0$.
The third equality is clear. The second one is obtained by counting the admissible $T$-covers over a curve
of genus $i-1$ and by multiplying for the $2$ possible gluing factors at the node. The first equality
is obtained by difference;
\item $A_i^{T,N}\cdot \delta_{j,c_3}^{S_3,S_3}=0$ for all $j$ and $A_i^{T,N}\cdot \delta_j'=0$ for all $j\neq i$
 and $A_i^{T,N}\cdot\delta_i'=-(2^{2i}-1)$ because it is the same of $A_i^{T,N}\cdot\pi^*\delta_i$.
\end{itemize}

Consider an effective divisor $ E$ of $\R_{g,S_3}^{S_3}$, as the pencils $A_i^{T,N}$ fill up the boundary divisor
$\Delta_i^{T,N}$ for $i\neq 10$, we have
$A_i^{T,N}\cdot \overline E\geq 0$ that is, by the relations above,
$$(i+1)\cdot a-\frac{2^{2i-1}+1}{2^{2i}-1}\cdot (6i+18)\cdot b_0'-\frac{2^{2i-1}-2}{2^{2i}-1}\cdot(6i+18)\cdot b_{0,c_2}+b_i'\geq 0,$$
which implies $b_i'\geq 3$ for all $i\geq 1$. The same is true for $i=10$ if
$\overline E$ does not contain the locus $\mc Z'$.\newline

We define $d:=\deg(A_{i,c_3}^{S_3,S_3})$.
As in the previous case, also for $A_{i,c_3}^{S_3,S_3}$ we have the equalities 
\begin{itemize}
\item $A_{i,c_3}^{S_3,S_3}\cdot \lambda=d\cdot (i+1)$;
\item $A_{i,c_3}^{S_3,S_3}\cdot \delta_{i,c_3}^{S_3,S_3}=-d$;
\item $A_{i,c_3}^{S_3,S_3}\cdot \delta_j'=0$ for all $j$
and $A_{i,c_3}^{S_3,S_3}\cdot \delta_{j,c_3}^{S_3,S_3}=0$ for all $j\neq i$. 
\end{itemize}
Similarly we have
$A_{i,c_3}^{S_3,S_3}\cdot (\delta_0'+\delta_{0,c_2}+\delta_{0,c_3})=A_{i,c_3}^{S_3,S_3}\cdot \pi^*\delta_0=d\cdot(6i+18)$.

The following inequality is also true
\begin{equation}\label{eq_00}
 A_{i,c_3}^{S_3,S_3}\cdot(\delta_{0,c_2}+\delta_{0,c_3})\geq A_{i,c_3}^{S_3,S_3}\cdot \delta_0'.
\end{equation}

For any effective divisor $E$ of $\R_{g,S_3}^{S_3}$, as the pencils $A_{i,c_3}^{S_3,S_3}$ fill 
up a connected component of the boundary divisor $\Delta_{i,c_3}^{S_3,S_3}$ for $i\neq 10$, we have
$A_{i,c_3}^{S_3,S_3}\cdot\overline E\geq 0$. This implies, by the inequalities above,
$$d\cdot b_{i,c_3}^{S_3,S_3}\geq A_{i,c_3}^{S_3,S_3}\cdot (2\delta_0'+3\delta_{0,c_2}+4\delta_{0,c_3})-d\cdot(i+1)\cdot a\geq $$
$$\frac{5}{2}\cdot A_{i,c_3}^{S_3,S_3}\cdot(\delta_0'+\delta_{0,c_2}+\delta_{0,c_3})-d\cdot(i+1)\cdot a> d\cdot 7.$$
The same is true if $i=10$ and $\overline E$ does not contain the locus $\mc Z'$.\newline

It remains to prove the inequality (\ref{eq_00}). Consider a $1$-marked curve $(C_1,p_*')$ of genus $i$ with
an autointersection node $q$, we call $\overline C_1$ the partial normalization at $q$.
By Remarks \ref{rmk_mono} and~\ref{mono2}, an admissible $S_3$-cover on $(C_1,p_*')$
is equivalent to an admissible $S_3$-cover on $(\overline C_1,p_*')$ plus the data of the $S_3$-type $\lq h_1\rq$
at $q$ and a gluing factor $h_{\gamma}$ at the same node (defined up to conjugation), such that $h_\gamma$ is in the centralizer of $h_1$.

Consider the pencil $A_{i,c_3}^{S_3,S_3}$ obtained by joining any such curve $(C_1,p_*')$ to a fixed
general $1$-marked curve $(C_2,p_*)$ of genus $g-i$ with a connected admissible $S_3$-cover on it. Any curve
in the (finite) intersection $A_{i,c_3}^{S_3,S_3}\cap \Delta_{0,c_2}$ induces a curve 
in the intersection $A_{i,c_3}^{S_3,S_3}\cap \Delta_0'$ by putting a trivial $S_3$-type at $q$ instead of $\lq h_1\rq$.
The same is true for
$A_{i,c_3}^{S_3,S_3}\cap\Delta_{0,c_3}$. All the points of $A_{i,c_3}^{S_3,S_3}\cap \Delta_0'$
are obtained at least once via these operations of $S_3$-type trivialization. This proves the inequality.\fine

\begin{prop}\label{prop_big}
If
$$K_{\overline\R_{g,S_3}}=a'\cdot \lambda + E'$$
on $\widetilde \R_{g,S_3}^{S_3}$, where $a'$ is a positive coefficient
and $E'$ an effective divisor not containing~$\mc Z'$,
then the canonical divisor
is big over the space $\overline\R_{g,S_3}^{S_3,\reg}$.
\end{prop}

\proof The equation in the hypothesis implies that
the canonical divisor is big on~$\widetilde\R_{g,S_3}^{S_3}$, because
$\lambda$ is a big divisor. If we consider the cloture
$\overline E'$ of $E'$ on the space~$\overline\R_{g,S_3}^{S_3,\reg}$, 
we observe that after Lemma \ref{lemma_kk}, it
respects the hypothesis of Lemma \ref{lemma_eff},
and therefore there exists another effective boundary divisor $E''\in\Pic_{\Q}(\overline\R_{g,S_3}^{S_3,\reg})$
such that
$K_{\overline\R_{g,S_3}}=a'\cdot\lambda+\overline E'+E''$
on $\overline\R_{g,S_3}^{S_3,\reg}$, and
the proof is completed.\fine

\section{Extension of pluricanonical forms}\label{ext}

We recall that to every twisted $G$-cover $(\ssC,\phi)$ is
uniquely associated an admissible $G$-cover $F\to C$ and \cor{vice versa}, see Remark \ref{rmk_equiv}.
In this and the following sections, we will mainly use 
 the notation of twisted $G$-covers.

In order to evaluate the Kodaira dimension of any moduli space $\overline\R_{g,G}$, 
we want to prove an extension result of pluricanonical forms, as done
for example by Harris and Mumford for $\overline\M_g$ (see \cite{harmum82})
and by Chiodo and Farkas for $\overline\R_{g,\ell}^0$ with $\ell<5$ and $\ell=6$ (see \cite{chiofar12}). In particular given a desingularization
$\widehat\R_{g,G}\to \overline\R_{g,G}$, and denoting by $\overline\R_{g,G}^{\reg}$ the sublocus of
regular points, we know that  $H^0\left(\widehat \R_{g,G},nK_{\widehat \R_{g,G}}\right)\subset H^0\left(\overline\R_{g,G}^{\reg},nK_{\overline\R_{g,G}}\right)$
and we would like to prove
$$H^0\left(\widehat \R_{g,G},nK_{\widehat \R_{g,G}}\right)=H^0\left(\overline\R_{g,G}^{\reg},nK_{\overline\R_{g,G}}\right)$$
for $n$ sufficiently big and divisible. This condition is verified locally
for smooth points and canonical singularities, it remains to treat the non-canonical locus
$\sing^{\nc}\overline\R_{g,G}$.\newline

First, we recall the structure of $\sing^{\nc}\overline\M_g$. Consider a curve $[C]\in\overline\M_g$,
then $[C]$ is a non-canonical singularity if and only if $C$ admits an elliptic tail automorphism of order $6$.
This means that $C$
has an irreducible component which is
an elliptic curve $E$ such that $E\cap \overline{C\backslash E}$ is a single point and $E$ admits
an automorphism of order $6$.
In the following, we denote by $(E,p_*)$ the elliptic tail of a curve $C$, where
$p_*$ is the preimage of the node. Moreover we define $C_1:=\overline{C\backslash E}$ and we mark the
preimage $p_*'$ of the node on $C_1$. With this notation $C=C_1\cup E$.\newline

Following again \cite{gale19}, we introduce the $T$-curves over $\R_{g,G}$.

\begin{defin}[$T$-curve]
A twisted $G$-cover $(\ssC,\phi)$ is a $T$-curve if
there exists
an automorphism $\msf a\in \underline\Aut(\ssC,\phi)$ such that its 
coarsening $a$ is an elliptic tail automorphism of order~$6$.
The locus of $T$-curves in $\overline\R_{g,G}$ is denoted by $T_{g,G}$.
\end{defin}

In order to introduce $J$-curves, we recall the definition of the age invariant associated to a
linear automorphism. Given a finite order automorphism $\msf h\in\GL(m)$,
we can diagonalize it as $\msf h=\diag\left(\xi_{r_1}^{a_1},\dots,\xi_{r_m}^{a_m}\right)$,
where $\xi_{r_i}$ is an $r_i$th root of the unit.
Then, its age is defined as $\age(\msf h):=\sum(a_i\slash r_i)$. A finite subgroup of $\GL(m)$ with
no quasireflections is junior if it contains a non-trivial element $\msf h$ with $\age(\msf h)<1$. While
the notion of age depends on the choice of the roots $\xi_{r_i}$, the notion
of junior group is independent from this choice. For a wider introduction to the age invariant see \cite[\S 5.1.1]{gale19}.

\begin{defin}[$J$-curve]
A twisted $G$-cover $(\ssC,\phi)$ is a $J$-curve if
the group 
$$\underline\Aut_C(\ssC,\phi)\slash\QR_C(\ssC,\phi),$$
which is the group of ghosts quotiented by its subgroup of quasireflections, is junior (see
\cite[Definition 4.1.2]{gale19} for the notion of ghost automorphism).
The locus of $J$-curves in $\overline\R_{g,G}$ is denoted by $J_{g,G}$.
\end{defin}

As proved in \cite[Theorem 5.1.8]{gale19}, $\sing^{\nc}\overline\R_{g,G})$ is the union of
$T_{g,G}$ and $J_{g,G}$. We
consider the cases where the second one is empty, and therefore we
treat the $T$-locus by generalizing the Harris-Mumford technique for~$\overline\M_g$.
In particular, we focus in the case of $G$ abelian group and in the case of $G=S_3$ the symmetric
group of order $3$.

\subsection{The case $G$ abelian group}

\begin{lemma}\label{lem_tcur}
Consider a twisted $G$-cover $(\ssC,\phi)$ with $G$ finite abelian group.
If $(\ssC,\phi)$ is a $T$-curve, then the restriction of the cover
to the elliptic tail is trivial.
\end{lemma}
\proof
By Proposition \ref{prop_gfg2}, the set of admissible $G$-covers over the
elliptic curve $(E,p_*)$ (that is the set of twisted $G$-covers over the same curve) is in bijection with the set of maps $\varpi\colon \pi_1(E_{\gen},p_*)\to G$.
If $(E,p_*)$ admits an order $6$ automorphism~$\msf a_6$, then
$E\cong \C\slash (\Z\oplus\Z\cdot\Omega)$,
where $\Omega$ is a primitive $6$th root of the unit
and $p_*$ is the origin. Therefore
$\msf a_6$ acts on $E$ as multiplication by $\Omega$.
The fundamental group $\pi_1(E,p_*)\subset \pi_1(E_{\gen},p_*)$ is generated by $a$ and $b$
which are the classes of the two laces $\gamma_a$ and $\gamma_b$ such that
\begin{align*}
\gamma_a\colon[0,1]\to \C &:\ \gamma_a(t)=t\\
\gamma_b\colon[0,1]\to \C &:\ \gamma_b(t)=t\cdot \Omega.
\end{align*}
We have as a consequence 
$\msf a_6(a)=b$, and $\msf a_6(b)=ba^{-1}$.
Therefore if we call $\varpi'$ the map $\msf a_6^*\varpi$,
by Proposition \ref{prop_gfg2}, 
$\msf a_6$ lifts to the cover 
if and only if
$\varpi'=\varpi$.
This is true if and only if 
$\varpi\equiv 1$,
\cor{i.e.}~the restriction of $(\ssC,\phi)$ to the elliptic
tail must be trivial.\fine

\begin{teo}\label{teoh0}
In the case of a moduli space $\overline\R_{g,G}$ of twisted $G$-covers
with $G$ finite abelian group,
we consider a desingularization $\widehat\R_{g,G}\to\overline\R_{g,G}$. If
the locus $J_{g,G}\subset\overline\R_{g,G}$ is empty,
then 
$$H^0(\overline\R_{g,G}^{\reg},nK_{\overline\R_{g,G}^{\reg}})=H^0(\widehat\R_{g,G},nK_{\widehat\R_{g,G}}),$$
for $n$ sufficiently big and divisibile.
\end{teo}

In \cite{harmum82} the same is proved for the moduli space $\overline\M_g$.
The idea is the following. Consider a general non-canonical singularity of $\overline\M_g$,
that is a point $[C]$ where $C=C_1\cup E$, $C_1$ is smooth of genus $g-1$ and without automorphisms,
and $E$
is an elliptic tail admitting an order $6$ automorphism.
Consider the operation of gluing any elliptic tail $E'$ at $C_1$ along the same node.
This gives an immersion
\begin{equation}\label{eq_psi}
\Psi\colon \overline\M_{1,1}\hookrightarrow\overline\M_g,
\end{equation}
and the image of $\Psi$ passes through the point $[C]$.
Furthermore, there exists a neighborhood $S=S([C])$ of $\im\Psi$ in $\overline\M_g$ with the following properties:
\begin{enumerate}
\item it exists a smooth $(3g-3)$-dimensional variety $B$ and
a birational morphism $g\colon S\to B$;
\item it exists a subvariety $Z\subset B$ of codimension $2$ such that
$g^{-1}(B\backslash Z)\cong B\backslash Z$;
\item as $B\backslash Z\subset S\subset \overline\M_g$, we have $B\backslash Z\subset\overline\M_g^0\subset\overline\M_g^{\reg}$,
where $\overline\M_g^0$ is the subspace of stable curves with trivial automorphism group.
\end{enumerate} 
$$
\begin{tikzcd}
\im\Psi \ar[r,hook] & S([C])\ar[d, "g"']\ar[rr,hook]&&\overline\M_g\\
&B & B\slash Z\ar[l,hook]\ar[ul, hook]\ar[r,hook]&\overline\M_g^0\ar[u,hook]
\end{tikzcd}
$$

This allows to conclude. Indeed, for every pluricanonical form $\omega$ on $S([C])^{\reg}$, we consider
its restriction to $B\backslash Z$, this extends to the smooth variety $B$ and pullbacks
to a desingularization $\widehat S([C])\to S([C])\to B$.
We use the Ludwig approach for $\overline\R_{g,\mmu_2}$ developed in \cite{lud10}.
 In order to complete this, we need a
generalization of the age tools and the age criterion \cite[Proposition 5.1.4]{gale19}.

\begin{prop}[see Appendix $1$ to $\S1$ of {\cite{harmum82}}]\label{genagec}
Consider a complex vector space $V\cong \C^n$, $\G\subset \GL(V)$ finite subgroup,
a desingularization $\widehat{V\slash\G}\to V\slash\G$ and a $\G$-invariant
pluricanonical form $\omega$ on $V$. Consider an element $\msf h$ in $\G$,
 $V^0\subset V$ the subset where $\G$ acts freely and
$\fix(\msf h)\subset V$ the fixed point set of $\msf h$. By an
abuse of notation we denote by $\fix(\msf h)$ also the image of the fixed
point set in $V\slash \G$. Let $U\subset V\slash \G$ be an open subset
such that $V^0\slash\G\subset U$ and such that for every $\msf h$ 
with $\age\msf h<1$ (with respect to some primitive root of the unit),
the intersection $U\cap\fix(\msf h)$ is non-empty.
We denote by $\widehat U\subset\widehat {V\slash\G}$ the preimage
of $U$ under the desingularization.

If $\omega$, as a meromorphic form on $\widehat{V\slash \G}$, is
holomorphic on $\widehat U$, then it is holomorphic on $\widehat{V\slash\G}$.
\end{prop}

\proof[Proof of Theorem \ref{teoh0}]
Consider a pluricanonical form $\omega$ on $\overline\R_{g,G}^{\reg}$.
We show that $\omega$ lifts to a desingularization of an open neighborhood
of every point $[\ssC,\phi]$ of $\overline\R_{g,G}$.

If $[\ssC,\phi]$ is a canonical singularity this is obvious by definition.

If $[\ssC,\phi]$ is a non-canonical singularity, at first we consider the case of a general non-canonical singularity.
As $\sing^{\nc}\overline\R_{g,G}=T_{g,G}$, then by Lemma \ref{lem_tcur}
a general point $[\ssC,\phi]\in T_{g,G}$ is a $T$-curve $\ssC$ whose coarse space
$C$ has two irreducible components $(E,p_*)$, an elliptic tail, and $(C_1,p_*')$ of genus $g-1$.
Moreover, if $(\ssC_1,\phi_1)$ and $(\msf E, \phi_{\msf E})$ are the restrictions,
then $(\msf E,\phi_{\msf E})$ is the trivial cover of $E$, $\ssC_1=C_1$ and 
$\underline\Aut(\ssC_1,\phi_1)$
is trivial.

Once we fix the twisted $G$-cover $(C_1,\phi_1)$, we consider the morphism 
$\Psi_1\colon \overline\M_{1,1}\to \overline\R_{g,G}$
sending any point $[E']$ of $\overline\M_{1,1}$ to the point $[C,\phi]$
obtained by joining $C_1$ and $E'$ along their marked points,
and by considering the $G$-cover $\phi$ such that $\phi|_{C_1}=\phi_1$ and $\phi|_{E'}$ is trivial.

Following \cite{lud10}, we see that the projection $\pi\colon\overline\R_{g,G}\to \overline\M_g$
sends $\im \Psi_1$ isomorphically on $\im \Psi$, and $\pi|_{\im\Psi}$ is a local isomorphism.
Indeed, $\defo(C;\phi)=\defo(C)$ and for every point of $\im\Psi_1$ the automorphism group $\underline\Aut(C,\phi)$
is isomorphic to $\Aut(C)$. Therefore if we
consider the neighborhood $S([C])$ of $[C]\in\im\Psi$ introduced by Harris and Mumford,
then up to shrinking $\pi^{-1}S([C])\cong S([C])$, we have a neighborhood of $[C,\phi]$
with the same properties.\newline

It remains to consider the case of any non-canonical singularity $[\ssC,\phi]$.
Here $\ssC$ is a twisted curve such that
$\ssC=\ssC_1\cup\bigcup\msf E^{(i)}$,
where the $\msf E^{(i)}$ are all the elliptic tails  of $\ssC$ admitting an elliptic tail automorphism of order~$6$. Again we follow
the last part of the Ludwig's demonstration of \cite[Theorem 4.1]{lud10}. We consider for each $i$
a small deformation $(\ssC^{(i)},\phi^{(i)})$ of $(\ssC,\phi)$ which fixes
the $i$th elliptic tail. That is,
$\ssC^{(i)}=\ssC_1^{(i)}\cup\msf E^{(i)}$ where $\ssC_1^{(i)}$ is irreducible. Moreover, the twisted $G$-cover admits no
non-trivial automorphism over $\ssC_1^{(i)}$ and it is unchanged over $\msf E^{(i)}$.
By the previous point we consider $S^{(i)}:=S([\ssC^{(i)},\phi^{(i)}])$. Up to shrinking the open
subsets $S^{(i)}$, they are all disjoint. 
Given the point $[\ssC,\phi]$ of $\overline\R_{g,G}$, we consider the local picture
of its universal deformation
$V:=\Def(\ssC,\phi)\cong \C^{3g-3}$, and recall that the local picture at $[\ssC,\phi]$ is
the same of
$V\slash\underline\Aut(\ssC,\phi)$
at the origin. We define
$$S([\ssC,\phi]):=(V\slash\underline\Aut(\ssC,\phi))\cup\left(\bigcup S^{(i)}\right).$$
If $V^0$ is the $V$ subset where $\underline\Aut(\ssC,\phi)$ acts freely,
let $U\subset S([\ssC,\phi])$ be the set
$$U:=(V^0\slash\underline\Aut(\ssC,\phi))\cup\left(\left(V\slash\underline\Aut(\ssC,\phi)\right)\cap\bigcup S^{(i)}\right).$$
If $\omega$ is a pluricanonical holomorphic form on $S^{\reg}$, then it extends to $\widehat S^{(i)}$ by definition of the neighborhoods $S^{(i)}$.
Moreover, by applying Proposition \ref{genagec} to the subset $U$, $\omega$ extends to $\widehat{\frac{V}{\underline\Aut(\ssC,\phi)}}$, and
therefore to the whole $\widehat S$.\fine

\begin{rmk}\label{rmk_roooo}
If we consider $G=\mmu_{\ell}$, we observe that this is the case (treated in \cite{chiofar12})
 of the moduli space $\overline\R_{g,\ell}^0$
of curves equipped  with a line bundle which is an $\ell$th root of the trivial bundle. The proof above applies
with minor adjustments also to the case of the moduli space $\overline\R_{g,\ell}^k$
for any $k$ (see \cite{gale15}), that is the moduli space of curves with an $\ell$th root of $\omega^k$. Therefore the extension result is true also in this case
if the $J$-locus $J_{g,\ell}^k$ of $\overline\R_{g,\ell}^k$ is empty.\newline
\end{rmk}

\subsection{The case $G=S_3$}

We know from \cite{gale19} that in $\overline\R_{g,S_3}$ the non-canonical singular locus
coincides with the $T$-locus.

\begin{lemma}\label{lem_Nm}
Consider a twisted $S_3$-cover $(\ssC,\phi)$ which
is a $T$-curve. 
If $\msf E$ is an elliptic tail admitting an elliptic tail automorphism $\msf a_6$ of order $6$,
then the restriction $(\msf E,\phi_{\msf E})$
 is a trivial cover
or an admissible $N$-cover, \cor{i.e.}~in this last case it has two connected components
and trivial $S_3$-type on the marked point $p_*$.
\end{lemma}
\proof We follow the same approach of Lemma \ref{lem_tcur}. By Proposition \ref{prop_gfg2},
the set of admissible $S_3$-covers on $E$ is in bijection with the conjugacy classes of maps
$\varpi\colon \pi_1(E_{\gen},p_*)\to S_3$. As before we note that if $(E,p_*)$ is an elliptic
curve admitting and order $6$ automorphism, then $E\cong \C\slash (\Z\oplus\Z\cdot \Omega)$, where $\Omega$ is
a primitive $6$th root of the unit and $p_*$ is the origin. Here $\msf a_6$ acts as multiplication by $\Omega$. 
If $a$ and $b$ are defined as before, then again $\msf a_6(a)=b$ and $\msf a_6(b)=ba^{-1}$.
If $\varpi'=\msf a_6^*\varpi$, then by Proposition \ref{prop_gfg2}, $\msf a_6$ lifts to the cover if and only if $\varpi'=h\cdot \varpi\cdot h^{-1}$ for
some element $h\in S_3$. This is true in the case $\varpi\equiv 1$ and for every
group morphism $\varpi\colon \pi_1(E)\twoheadrightarrow N\subset S_3$, where $N$ is the normal non-trivial
subgroup of $S_3$, that is the group generated by any $3$-cycle. By definition of admissible $N$-cover, the thesis follows.\fine

As observed in \S\ref{sec_RRR}, there exists an isomorphism, $\overline\R_{1,S_3}^1\cong \overline\M_{1,1}$. 
We remarked in the same section that there exists a natural $2:1$ morphism $\overline\R_{1,\mmu_3}^{\mmu_3}\to\overline\R_{1,S_3}^N$.
As the natural morphism $\overline\R_{1,\mmu_3}^{\mmu_3}\to\overline\M_{1,1}$ is $8:1$, we have that $\Psi_{N}\colon\overline\R_{1,S_3}^N\to \overline\M_{1,1}$
is a $4:1$ morphism.

\begin{rmk}
Consider $(E,p_*)$ a general elliptic curve, we list the four
preimages of $[E,p_*]\in\overline\M_{1,1}$ with respect to $\Psi_N$.
We characterize every class with a representative $\varpi\colon\pi_1(E,p_*)\twoheadrightarrow N$ of the correspondent\
conjugacy class of morphisms. We recall that $\pi_1(E,p_*)$ is a free abelian group generated by $a$ and $b$.
\begin{center}
{\renewcommand\arraystretch{1.4}

\begin{tabular}{|c|c|}

& $\varpi\colon\pi_1(E,p_*)\twoheadrightarrow N$  \\ \hline
 (i) &$a\mapsto1;\ \ b\mapsto (123)$ \\ \hline
(ii) &$a\mapsto(123);\ \ b\mapsto(123)$ \\ \hline
(iii) & $a\mapsto(123);\ \ b\mapsto(132)$ \\ \hline
(iv) & $a\mapsto(123);\ \ b\mapsto 1$ \\ \hline
\end{tabular}}
\end{center}
\end{rmk}

We recall that $\overline\M_{1,1}\cong\p^1$, $\overline\R_{1,S_3}^N$ is a connected curve and we are interested
in finding the branch points for $\Psi_N$. 
We observe that the only automorphism of a general
elliptic curve is the natural involution $i$, and this always lifts to any admissible $S_3$-cover. 
Following Remark~\ref{rmk_local}, we detect a branch point in three cases.
\begin{itemize}
\item The elliptic curve $E_4$ admitting an automorphism $\msf a_4$ of order $4$. Automorphism $\msf a_4$
acts on $\pi_1(E_4,p_*)$ by sending $a\mapsto b$ and $b\mapsto a^{-1}$, therefore it does not lift
to any admissible $N$-cover: in particular it exchanges classes $(i)$ and $(iv)$, and also classes $(ii)$ and $(iii)$.
Therefore over $[E_4]\in \overline\M_{1,1}$ we have two branch points of order $2$ that we denote by $[E_4']$ and $[E_4'']$.
\item The elliptic curve $E_6$ admitting an automorphism $\msf a_6$ of order $6$. Automorphism $\msf a_6$ acts
on $\pi_1(E_6,p_*)$ by sending $a\mapsto b$ and $b\mapsto ba^{-1}$, therefore it lifts to the cover $(iii)$. Moreover,
it sends $(i)$ to $(ii)$, $(ii)$ to $(iv)$ and $(iv)$ to $(i)$. Therefore over $[E_6]\in\overline\M_{1,1}$ we have one branch point $[E_6']$
of order $3$, and moreover we have another preimage point $[E_6'']$.
\item The curve $E_0$ with an autointersection node $q_1$, whose normalization $\overline E_0$ is a rational curve.
By Remark \ref{rmk_local}, $[E_0]\in\overline\M_{1,1}$ has two preimages via $\Psi_N$,
one is obtained by putting an order $3$ stabilizer at $q_1$, we call this preimage $[E_0']$
and it is an order $3$ branch point. The other preimage, called $[E_0'']$, is associated
to an admissible $N$-cover with trivial stabilizer at $q_1$.
\end{itemize}

\begin{lemma}\label{lem_rat}
The moduli space $\overline\R_{1,S_3}^N$ is isomorphic to $\p^1$.
\end{lemma}
\proof This is a consequence of Hurwitz formula plus the observations
we just listed.\fine

We consider a general genus $g-1$ curve with one marked point $(C_1,p_*')$ and
a twisted $S_3$-cover $(C_1,\phi_1)$ with trivial $S_3$-type at $p_*'$.
For any twisted $G$-cover $(\msf E,\phi_{\msf E})$ which is an admissible $N$-cover over an elliptic curve $(E,p_*)$, 
we glue $(C_1,\phi_1)$ and $(\msf E,\phi_{\msf E})$ along their marked points, to get a node $q$. We obtain a map
$\Upsilon\colon\overline\R_{1,S_3}^N\to \overline\R_{g,S_3}$
which is an isomorphism into the image and fits into the diagram
$$
\begin{tikzcd}
\overline\R_{1,S_3}^N\ar[d,"\Psi_N"']\ar[r,hook, "\Upsilon"] & \overline\R_{g,S_3}\ar[d,"\pi"]\\
\overline\M_{1,1}\ar[r,hook,"\Psi"'] &\overline\M_g.
\end{tikzcd}
$$

Therefore the projection $\pi|_{\im\Upsilon}\colon\im\Upsilon\to \im\Psi$ is the $4:1$ morphism
we described above. 

\begin{rmk}\label{rmk_defdef}
Before stating the main theorem of this section, we describe the local
picture of $\pi$ at the points of $\im\Upsilon$.
We call $(\ssC,\phi)$ the admissibile $S_3$-cover obtained via the gluing.
By Remark \ref{rmk_local}, the local picture of $\overline\R_{g,S_3}$ at 
$[\ssC,\phi]$ is
$$\left(\defo(C_1,p_*')\oplus\A^1_{\tilde t_q}\oplus\A^1_{\tilde t_E}\right)\slash\underline\Aut(\ssC,\phi).$$
Here $\tilde t_q$ and $\tilde t_E$ are the (non-canonical) coordinates associated to the smoothing of node $q$ and
to the deformation of $(\msf E,\phi_{\msf E})$. At the same time the local picture of $\overline\M_g$ at $[C]$ is
$$\left(\defo(C_1,p_*')\oplus\A^1_{t_q}\oplus\A^1_{t_E}\right)\slash\Aut(C).$$
Here $t_q$ and $t_E$ are associated to the smoothing of node $q$ and the deformation of the elliptic curve $E$.
As the $S_3$-type at $q$ is trivial, $\tilde t_q=t_q$, while for the other coordinate there are different cases.
As the automorphism group acts non-trivially only on $t_q, t_E$ and $\tilde t_E$, we focus on these coordinates.
Moreover we recall that the canonical elliptic tail involution of $E$, acts trivially.
\end{rmk}

We consider again the case treated by Harris and Mumford in \cite{harmum82}, the birational morphism $g\colon S([C])\to B$
at a general non-canonical singularity $[C]$. This means that the associated curve $C$ is the junction
of $(C_1,p_*')$, a genus $g-1$ smooth automorphism free curve, and $(E,p_*)$, an elliptic tail.
By construction $S([C])\cong S_1\times \defo(C_1,p_*')$ and $B\cong B_1\times \defo(C_1,p_*')$,
where $S_1$ is a singular surface, $B_1$ is a smooth surface and $g$ is the identity
on the component $\defo(C_1,p_*')$. We denote by $D$ the projection of $\im\Psi$ on $S_1$.
In particular the coordinates $t_q,t_E$ span exactly the tangent space to $S_1$ at any point
 of $D$.

If $[\ssC,\phi]$ is a general non-canononical singularity or $\overline\R_{g,S_3}$, then
the curve coarse space $C$ is  a union $C_1\cup E$ as before. We denote by $S([\ssC,\phi])$
the component of $\pi^{-1}(S([C]))$ containing $[\ssC,\phi]$. By Lemma \ref{lem_Nm},
 the restriction $\phi_{\msf E}$
is trivial or an admissible $N$-cover, we focus in the second case and
we observe that as a consequence $\im\Upsilon\subset S([\ssC,\phi])$.
Following the construction above, $S([\ssC,\phi])=S_2\times\defo(C_1,p_*')$,
where $S_2$ is a singular surface. We call $D'$ the projection of $\im\Upsilon$ on $S_2$.
With an abuse of notation we call $\pi\colon S_2\to S_1$ the natural projection.

In the following list we describe the local picture of $S_2\to S_1$ at any point.
For the sake of simplicity we identify $D'$ with $\overline\R_{1,S_3}^N$
and $D$ with $\overline\M_{1,1}$.
We use the notation $\left(\frac{a_1}{r_1},\cdots,\frac{a_m}{r_m}\right)$ to denote the diagonal matrix $\diag(\xi_{r_1}^{a_1},\dots,\xi_{r_m}^{a_m})$,
where $\xi_r=\exp(2\pi i\slash r)$ is the privileged $r$th root of the unit. 

\begin{itemize}
\item At the two points $[E_4'],[E_4'']\in S_2$, the automorphism $\msf a_4$ does not lift.
It acts as $t_q\mapsto \xi_2t_q$ and $t_E\mapsto \xi_2t_E$ locally at $[E_4]\in S_1$, therefore the local picture of $\pi$ at $[E_4']$
and $[R_4'']$ is the canonical projection
$\C^2\to \C^2\slash\left(\frac{1}{2},\frac{1}{2}\right)$.
\item At the two points $[E_6'],[E_6'']\in S_2$, the automorphism $\msf a_6$ does not lift to $[E_6']$ but it lifts
to $[E_6'']$. It acts as $t_q\mapsto \xi_3t_q$ and $t_E\mapsto \xi_3t_E$ at $[E_6]\in S_1$, therefore the local picture
of $\pi$ at $[E_6']$ is the projection
$\C^2\to \C^2\slash\left(\frac{1}{3},\frac{1}{3}\right)$,
and the local picture of $\pi$ at $[E_6'']$ is the identity
$\C^2\slash\left(\frac{1}{3},\frac{1}{3}\right)\to\C^2\slash \left(\frac{1}{3},\frac{1}{3}\right)$.
\item At the two points $[E_0'],[E_0'']\in S_2$, we have that $t_E=\tilde t_E^3$ in the case of $[E_0']$
because of
the definition of admissible $S_3$-cover. Instead in the case of $[E_0'']$, $t_E=\tilde t_E$. In both cases there
are no additional automorphisms, therefore the local picture of $\pi$ at $[E_0']$ is the projection
$\C^2\to \C^2\slash\left(1,\frac{1}{3}\right)\cong \C^2$,
and the local picture at $[E_0'']$ is the identity $\C^2\to \C^2$.
\item Elsewhere on $D'$ and on the whole surface $S_2$, the projection $\pi$ is étale, therefore its local
picture is the identity $\C^2\to\C^2$.
\end{itemize}

\begin{teo}\label{teo_extpcf}
Consider a desingularization $\widehat\R_{g,S_3}\to\overline\R_{g,S_3}$
of the moduli space of genus $g$ curves equipped with an admissible $S_3$-cover.
Then,
$$H^0\left(\overline\R^{\reg}_{g,S_3},nK_{\overline\R_{g,S_3}^{\reg}}\right)=H^0\left(\widehat\R_{g,S_3},nK_{\widehat\R_{g,S_3}}\right)$$
for $n$ sufficiently big and divisible.
\end{teo}
\proof We are going
to prove the result for a general non-canonical singularity of $\overline\R_{g,S_3}$. 
As a consequence, the extension of pluricanonical forms is true also for any non-canonical singularity
via a patchwork of ``good'' neighborhoods analogous to what we do in the proof of Theorem~\ref{teoh0}.

Consider a general non-canonical singularity $[\ssC,\phi]\in\overline\R_{g,S_3}$, therefore
$\ssC=C_1\cup \msf E$ as described before and by Lemma \ref{lem_Nm}
the restriction $\phi_{\msf E}$ is trivial or it is an admissible $N$-cover.
In the first case, $[\ssC,\phi]$ is in a component of $\pi^{-1}\left(\im\Psi\right)$ which is isomorphic
to $\im\Psi$ and the construction of $S([\ssC,\phi])$ is the same of Theorem \ref{teoh0}. In the second case,
$[\ssC,\phi]$ is in $\im\Upsilon$ and we are going to show a neighborhood $S([\ssC,\phi])=S_2\times \defo(C_1,p_*')$ and a birational morphism
$S_2\to B_2$ such that $B_2$ has at most canonical singularities.

Consider the blowup at the points $[E_4],[E_6]\in D\subset S_1$. We denote by $A_4$ and $A_6$ the two corresponding exceptional
curves, by the description
of singularities we gave above, the autointersection numbers of these divisors are $A_4^2=-2$ and $A_6^2=-3$.
We denote by $S_1^*$ the blown up surface and by $\overline D$ the strict transform of $D$ after the blowups. We know by the Harris-Mumford result
resumed above, that there exists a contraction of the curve $A_4\cup A_6\cup\overline D$, and the contracted surface
is smooth. By \cite[Theorem II.11]{beau96} we must have a sequence of $(-1)$-curve contractions, and this happens
if and only if $\overline D^2=-1$, and all the three curves are rational. 

Consider the blowup at the points $[E_4'],[E_4''],[E_6'],[E_6'']\in D'\subset S_2$.
We denote by $A_4',A_4'',A_6',A_6''$ the corresponding exceptional divisors,
by $S^*_2$ the blown up surface and
 by $\overline D'$ the strict transform of $D'$.
From the description of the singularities we gave above, we know $(A_4')^2=(A_4'')^2=(A_6')^2=-1$
and $(A_6'')^2=-3$. Moreover, after the blowups there exists a morphism $\tilde \pi\colon S^*2\to S^*1$
of degree $4$ and whose ramification locus is $2A_6'$ (plus some component of codimension $2$).
Therefore $\tilde\pi^*\overline D=4\overline D'$, and so $(\overline D')^2=-4$. By costruction all the $A_i'$ and $A_i''$
are rational curves and $\overline D'$ is a rational curve too by Lemma \ref{lem_rat}. We obtain that there exists a birational morphism
$S_2\to B_2$
contracting the curve $A_4'\cup A_4''\cup A_6'\cup A_6''\cup D'$. Indeed, by \cite[Theorem III.5.1]{bhpv04}
there exists such a contraction and $B_2$ 
 has a singularity of type $\C^2\slash\left(\frac{1}{2},\frac{1}{2}\right)$, that is a canonical singularity.

We obtained a contraction $g'\colon S([\ssC,\phi])\to B'=B_2\times \defo(C_1,p_*')$ such that $B'$ has only canonical singularities, moreover
there exists a locus $Z'\subset B'$ of codimension $2$ such that $(g')^{-1}(B'\backslash Z')\cong B'\backslash Z'$.
By construction $(g')^{-1}(B'\backslash Z')$ is naturally a subset of $\overline\R_{g,S_3}$ and this allows to conclude
as in the previous case. \fine

\section{Evaluating the Kodaira dimension}\label{finally}
In order to calculate the Kodaira dimension, we need to develop some calculations in the tautological ring of the moduli space.
 In \S\ref{sec_grr}
we develop Grothendieck Riemann-Roch type calculations for vector bundles.
In \S\ref{sec_big} we apply it to evaluate the canonical divisor $K_{\overline\R_{g,S_3}}$ and
prove its bigness over $\overline\R_{g,S_3}^{S_3}$, the connected
component of $\overline\R_{g,S_3}$ of  connected $S_3$-covers.

\subsection{Adapted Grothendieck Riemann-Roch}\label{sec_grr}

\subsubsection{Tautological classes}
In this section we recall some well known \cor{tautological} classes in the Chow ring $A^*(\overline\M_{g,n})$
of the moduli space of curves, and their generalizations to $\overline\R_{g,G}$.
For a wider survey of the tautological relations and the tautological rings structure
see~\cite[\S17]{acgh11}, \cite{pix13} and \cite{fabpan13}.

There exists two natural morphisms ``coming from the geometry of curves'' on the moduli spaces $\overline\M_{g,n}$.
With this we mean that we can define them using the
modular interpretation of the space.
\begin{itemize}
\item The \cor{forgetful morphism} is a morphism
$$\mu\colon\overline\M_{g,n}\to\overline\M_{g,n-1}$$
sending any geometric point $[C;p_1,\dots,p_n]$ to the same
marked stable curve without the last point, $[C;p_1,\dots,p_{n-1}]$.
\item The \cor{gluing morphisms} are of two types
$$\iota\colon\overline\M_{g_1,n_1+1}\times\overline\M_{g_2,n_2+1}\to \overline\M_{g_1+g_2,n_1+n_2}\ \ 
\m{and}\ \ 
\iota\colon\overline\M_{g-1,n+2}\to\overline\M_{g,n}.$$
In the first case $[C;p_1,\dots,p_{n_1+1}]\times[C';p_1',\dots,p_{n_2+1}']$
is sent two the 
junction of $C$ and $C'$ along the marked points $p_{n_1+1}$ and $p_{n_2+1}'$.
The new curve maintains all the other marked points. In the second
case $[C;p_1,\dots,p_{n+2}]$ is sent to the quotient curve $C\slash(p_{n+1}\sim p_{n+2})$
 with the same other marked points.\newline
\end{itemize}

\begin{defin}
The system of tautological rings $R^*(\overline\M_{g,n})\subset A^*(\overline\M_{g,n})$
with $g,n$ varying on the non-negative integers,
is the smallest system of $\Q$-algebras closed
under the pushforwards of the forgetful morphisms and the gluing morphisms.
\end{defin}

We define $n$ tautological $\psi$-classes inside the Chow ring $A^*(\overline\M_{g,n})$.
We consider the universal family
$u\colon\mc C_{g,n}\to \overline\M_{g,n}$,
where $\mc C_{g,n}$ is a Deligne-Mumford stack such that every geometric fiber of $u$ is isomorphic
to the associated $n$-marked stable curve, and there exist $n$ sections 
$\sigma_1,\dots,\sigma_n\colon\overline\M_{g,n}\to\mc C_{g,n}$.
For every $i=1,\dots, n$ the line bundle $\mc L_i$ over $\overline\M_{g,n}$ is the $i$th
cotangent line bundle, $\mc L_i:=\sigma_i^*\left(T_u^{\vee}\right)$.
Then we define
$$\psi_i:=c_1(\mc L_i)\in A^1(\overline\M_{g,n}).$$

\begin{rmk}
The $\psi$-classes are in the tautological ring $R^1(\overline\M_{g,n})$,
as showed for example in \cite{fabpan13}.
\end{rmk}

There are two other type of classes that are very important for our analysis and belong to the tautological ring.
To introduce $\kappa$-classes we consider the log-canonical line bundle on~$\mc C_{g,n}$,
$\omega_u^{\log}:=\omega_u(\sigma_1+\dots+\sigma_n)$.
Therefore we have
$$\kappa_d:=u_*\left(c_1(\omega_u^{\log})^{d+1}\right).$$
As before the class $\kappa_d$ is well defined in the Chow ring $A^d(\overline\M_{g,n})$.
We state without proof the known fact that
$\kappa_d=\mu_*(\psi_{n+1}^{d+1}$) 
and being $\psi_{n+1}$ in the tautological ring, the $\kappa$-classes
is contained in the tautological ring too.

The Hodge bundle over $\overline\M_{g,n}$ is the rank $g$ vector bundle
$\mc E:=u_*\omega_u$,
\cor{i.e.}~the vector bundle whose fiber at $[C;p_1,\dots,p_n]$ is $H^0(C;\omega_C)$.
The \cor{Hodge class}~is 
$$\lambda:=c_1(\mc E)\in A^1(\overline\M_{g,n}).$$
The Hodge class is proved to be in the tautological ring $R^1(\overline\M_{g,n})$ in \cite{mum83}.\newline

We define the universal family also on $\overline\R_{g,G}$ and we use the
same notation $u\colon \mc C_{g,G}\to \overline\R_{g,G}$. We will see in the next
section that this family is equipped with a universal twisted $G$-cover.
On $A^*\left(\overline\R_{g,G}\right)$ it is possible to define as before the $\kappa$ classes and the Hodge class $\lambda$.\newline

\subsubsection{Using Grothendieck Riemann-Roch in the Chow ring of $\mc R_{g,G}$}
The universal family of curves over $\mc C_{g,G}\to\overline\R_{g,G}$
is equipped with a universal twisted $G$-cover 
$\Phi\colon \mc C_{g,G}\to BG$.
In particular for every geometric point $[\ssC,\phi]$ of $\overline\R_{g,G}$,
the restriction of $\Phi$ to the associated geometric fiber
is isomorphic to the twisted $G$-cover $\phi\colon \ssC\to BG$.

We consider the singular locus $\mc N\subset\mc C_{g,G}$ of the universal family,
whose points are the nodes of $\mc C_{g,G}$ fibers.
Furthermore, we consider the stack $\mc N'$ 
whose points are nodes equipped with the choice of a privileged branch.
There exists
a natural étale double cover $\mc N'\to \mc N$, and an involution $\varepsilon\colon \mc N'\to \mc N'$ associated
to the branch switch.
There exists also a natural decomposition of $\mc N'$:
given a conjugacy class $\lq h\rq $ in~$\lq G\rq $, we denote by 
$\mc N_{i,\lq h\rq }'\subset \mc N'$ the substack
of nodes such that 
the associated privileged branch is in a component of genus $i$,
and it has $\lq h\rq $ as $G$-type. 
In the case of a node of type $0$, the component is
$\mc N'_{0,\lq h\rq }$.
Therefore
$$\mc N'=\bigsqcup_{
\substack{0\leq i\leq g-1,\\ 
\lq h\rq \in\lq G\rq } }\mc N_{i,\lq h\rq }'.$$
We denote the natural projection by $j\colon\mc N'\to\overline\R_{g,G}$. 
Furthermore, we denote by $j_{i,\lq h\rq }$ the restriction 
of the map $j$ to the component $\mc N'_{i,\lq h\rq }$.
We finally define the classes $\psi$ and $\psi'$
on~$\mc N'$:
$$\psi:=c_1(T_{\mc N'}^{\vee});\ \ \ \psi':=c_1(\varepsilon^*T_{\mc N'}^{\vee}).$$

In order to evaluate the Chern character
of some line bundle pushforwards, we state a generalization of the Grothendieck Riemann-Roch formula,
following the approach of Chiodo in \cite{chio08}. In the case treated by Chiodo, this vector bundle
is the universal root of the trivial (or canonical) line bundle. Here we generalize
to the case of any vector bundle coming from a representation of the group $G$.

Consider $W$, a dimension $w$ representation of group $G$, then $W$
can be regarded as a vector bundle on $BG$. We consider
the universal cover $\Phi\colon \mc C_{g,G}\to BG$, then the
pullback 
$$W_{\mc C_{g,G}}:=\Phi^*W$$
is a vector bundle of rank $w$ on the universal family $\mc C_{g,G}$.
From now on we will use the more compact notation $W_{\mc C}$.

We observe that the projection $j_{i,\lq h\rq }\colon \mc N'_{i,\lq h\rq }\to \overline\R_{g,G}$
is locally isomorphic to $B\mmu_r\to \Spec\C$
at every point of $\overline\R_{g,G}$, where $r=r(h)$
is the order of class $\lq h\rq $ in~$G$.
We follow the approach of \cite[\S2.2]{tse10} to decompose the restriction of $W_{\mc C}$
to any locus~$\mc N_{i,\lq h\rq }'$. 
The local picture of $\mc C_{g,G}$ at any point of $\mc N_{i,\lq h\rq }'$ is a Deligne-Mumford stack
$[U\slash \mmu_r]$ where $U$ is an affine scheme,
therefore the vector bundle restriction $W_{\mc C}|_{\mc N_{i,\lq h\rq }'}$
is a $\mmu_r$-equivariant vector bundle. Following \cite[\S2.3.2]{gale19}, $\mmu_r$ acts naturally on $W_{\mc C}|_{\mc N_{i,\lq h\rq }'}$,
so there exists a decomposition in subbundles,
$$\left.W_{\mc C}\right|_{\mc N_{i,\lq h\rq }'}=W_0\oplus W_1\oplus\cdots\oplus W_{r-1},$$
where $W_k$ is the eigen-subbundle with eigenvalue $\xi_r^k$ for all $k=0,\dots,r-1$.
We denote by $w_{i,\lq h\rq }(k)$ the rank of $W_k$, and clearly these integers
satisfy the equation 
$\sum_{k}w_{i,\lq h\rq }(k)=w$.
We also recall the Bernoulli polynomials $B_d(x)$ defined by the 
generating function
$$\frac{te^{xt}}{e^t-1}=\sum_{d=0}^{\infty} B_d(x)\frac{t^d}{d!}.$$
The Bernoulli numbers $B_d:=B_d(0)$ are the evaluations
of the Bernoulli polynomials at $0$. 
With this setting we can state the following.

\begin{prop}\label{prop_chern}
On $\overline\R_{g,G}$ we have the following
evaluation for the degree $d$ component
of the Chern character of $ R u_*W_{\mc C}$.

$$\ch_d(Ru_*W_{\mc C})=\frac{w\cdot B_{d+1}}{(d+1)!}\kappa_d+$$
$$+\frac{1}{2}\cdot \sum_{\substack{0\leq i\leq g-1,\\ 
\lq h\rq \in\lq G\rq } } \left(\sum_{0\leq k< r(h)} \left(\frac{r(h)^2\cdot w_{i,\lq h\rq }(k)\cdot B_{d+1}\left(k\slash r(h)\right)}{(d+1)!}\right)\cdot(j_{i,\lq h\rq })_*\left(\sum_{a+a'=d-1}\psi^a(-\psi')^{a'}\right)\right).
$$
\end{prop}

This formula follows directly from Tseng formula (7.3.6.1) in \cite{tse10}.
In Tseng notation the morphism $ev_{n+1}$ is the $u$ morphism,
the representation $W$ is denoted by $F$ and moreover
$$(u_*(\ch(ev^*W)Td^\we(L_{n+1})))_d$$
is the first term in our formula, the one with $\kappa$ classes. 
The $\psi_i$-classes terms are associated to marked points and therefore
are absent in our formula. Finally, the terms $A_m$ (see \cite[Definition 4.1.2]{tse10})
give the last term of our formula.\newline

\subsection{Bigness of the canonical divisor}\label{sec_big}

We consider the moduli space $\overline\R_{g,S_3}^{S_3}$
of genus $g$ curves equipped with a connected admissible $S_3$-cover.
This is a component of~$\overline\R_{g,S_3}$.

The goal of this section is to prove that the canonical divisor of $\overline\R_{g,S_3}^{S_3}$
is big over the subspace $\widetilde\R_{g,S_3}^{S_3}\subset\overline\R_{g,S_3}^{S_3}$
of covers on
 irreducible curves,
 for every
odd genus $g>11$. By Proposition \ref{prop_big}, this implies that
$\overline\R_{g,S_3}^{S_3}$ is of general type for every odd genus $g\geq 13$.

The approach follows the strategy of \cite{cefs13} for $\overline\R_{g,\mmu_3}^{\mmu_3}$.
We write down the canonical divisor as a sum
\begin{equation} \label{eq_K}
K=\alpha\cdot \U+\beta\cdot \mm+E+\gamma\cdot \lambda\in\Pic_\Q(\widetilde \R_{g,S_3}^{S_3}).
\end{equation}
Here $\U,\mm,E$ are effective divisors, $\lambda$ is the hodge class,
$\alpha,\beta$ positive coefficients and $\gamma$ a strictly positive coefficient. 

\subsubsection{Basic notions of syzygy theory}
The divisor $\U$ in equation (\ref{eq_K}) 
will be defined following the approach of Chiodo-Eisenbud-Farkas-Schreyer
paper \cite{cefs13}. It is the locus of curves
with ``extra'' syzygies with respect to a particular vector bundle,
and it has a determinantal structure over an open subset of
$\overline\R_{g,S_3}^{S_3}$.
To properly define this, we recall some fundamental notions
of syzygy theory over stable curves, following
the notations of Aprodu-Farkas paper~\cite{aprofar11}.

Consider a finitely generated graded module $N$ over the polynomial ring
$S=\C[x_0,\dots,x_n]$. The module has a minimal free resolution
$0\leftarrow N\leftarrow F_0\leftarrow F_1\leftarrow\cdots$,
where
$$F_i=\sum_j S(-i-j)^{b_{i,j}}.$$
The numbers $b_{i,j}$ are well defined
and are called the Betti numbers of $N$,
moreover we have
$$b_{i,j}=\dim (\Tor_i^S(N;\C))_{i+j}.$$

\begin{rmk}\label{rmk_nat}
In an irreducible flat family of modules $N_t$,
the Betti numbers $b_{i,j}(N_t)$ are semicontinuous,
and simultaneously take minimal values on an open set. 
The jumping locus for their values is where one of the values
$b_{i,j}(N_t)$ is bigger than this minimum.
\end{rmk}

For every stable curve we consider a line bundle $L\in \Pic(C)$,
a sheaf $\mc F$ on $C$, the polynomial ring $S:=\Sym H^0(C;L)$ and
 the graded $S$-module
$$\Gamma_C(\mc F;L):=\bigoplus_{n\in \Z}H^0(C;\mc F\xx L^{\xx n}).$$
From now on we will use the Green notation by calling
$$K_{i,j}(C;\mc F,L):=(\Tor_i^S(\Gamma_C(\mc F,L);\C))_{i+j}.$$
One fundamental idea of syzygy theory is that the vector space
$K_{i,j}(C;\mc F,L)$ can be evaluated
via a minimal $S$ resolution of $\C$, where the latter is seen
as a graded $S$-module. To do this we introduce the vector bundle
$M_L$ of rank $H^0(L)-1$ over~$C$, via the following
short exact sequence
\begin{equation}\label{eq_M}
1\to M_L\to H^0(C;L)\xx \Oo_C\xrightarrow{\ev} L\to 1.
\end{equation}
Using this vector bundle we can state the following lemma

\begin{lemma}[see {\cite[Theorem 2.6]{aprofar11}}]\label{cane} 
Consider $C$ a stable curve, $L$ a line bundle on it, $\mc F$ a coherent sheaf  on it and
$m$ a positive integer. Then,
$$K_{m,1}(C;\mc F,L)=H^0\left(\bigwedge^mM_L\xx \mc F\xx L\right).$$\newline
\end{lemma}

\subsubsection{The Koszul divisor}
The symmetric group $S_3$ has $3$ irreducible representations.
\begin{enumerate}
\item The trivial representation $\I\colon S_3\to \GL(1;\C)=\C^*$;
\item the parity representation $\epsilon\colon S_3\to \C^*$, which sends
even elements to $1$ and odd elements to $-1$;
\item given a vector space $R\cong \C^2$, the representation 
$\rho\colon S_3\to \GL(R)\cong \GL(2;\C)$
such that
$$\rho((12))=\left(\begin{matrix}
-1 & -1\\
0 & 1\end{matrix}\right)\ \ \ 
\rho((23))=\left(\begin{matrix}
0 &1\\
1& 0
\end{matrix}\right)\ \ \ 
\rho((13))=\left(\begin{matrix}
1 & 0\\
-1 & -1
\end{matrix}\right).
$$
\end{enumerate}

\begin{rmk}
In particular, if we consider the tautological representation over the vector
space 
$P:=\langle v_1,v_2,v_3\rangle_\C$
such that $S_3$ acts naturally by permutation, then $P$ is the direct sum of the 
trivial representation and $R$,
$P=\C_\I\oplus R$.
At the same time, if we consider the regular representation $\C[S_3]$ of dimension $|S_3|=6$,
we have
$\C[S_3]=\C_\I\oplus\C_\epsilon \oplus R^{\oplus 2}$.
\end{rmk}

As explained above, we consider $R$ as a line bundle over the stack $BS_3$.
Given a twisted $G$-cover $[\ssC,\phi]$, the pullback $R_{\ssC}:=\phi^*R$ is a 
rank $2$ vector bundle over $\ssC$. In particular if $\ssC=C$ is a scheme theoretic
curve, $R_C$ is a usual scheme theoretic vector bundle.\newline

We consider the case of odd genus $g=2i+1$, and
we focus on the Koszul cohomology $K_{i,1}$.
We introduce the locus $\mc U_g$
as the locus with non-zero cohomology~$K_{i,1}$.
Supposing that the minimal value assumed by $\dim K_{i,1}$
on $\overline\R_{g,S_3}^{S_3}$ is $0$,
this is therefore the jumping locus for $K_{i,1}$. In particular
we will show that $\mc U_g$ is a virtual divisor, \cor{i.e.}~$\mc U_g$
is an effective divisor if and only if $\dim K_{i,1}$ takes
value $0$ on a general curve, or equivalently on at least one curve.

\begin{defin}\label{def_U}
Given an odd genus $g=2i+1$,
$$\U_g:=\left\{[C,\phi]\in\R_{g,S_3}^{S_3}\ |\ K_{i,1}(C;R_C,\omega_C)\neq 0\right\}\subset \overline\R_{g,S_3}^{S_3}.$$
\end{defin}

We want to show that $\U_g$ is a virtual divisor.
For any stable curve $C$,
by Lemma \ref{cane} applied in the case $L=\omega_C$ with the vector bundle $R_C$ as the sheaf $\mc F$,
we have 
$$K_{i,1}(C;R_C,\omega_C)=H^0\left(\bigwedge^i M_{\omega}\xx \omega_C\xx R_C\right).$$
We can reformulate the definition of $\mc U_g$ with another
scheme theoretic condition.
Consider the Equation (\ref{eq_M}) in the case of the canonical line bundle $L=\omega_C$.
As $\omega_C$ is a line bundle, by a well known property we have
the short exact sequence
$$0\to \bigwedge^iM_{\omega}\xx \omega_C\xx R_C\to \bigwedge^iH^0(\omega_C)\xx \omega_C\xx R_C\to \bigwedge^{i-1}M_{\omega}\xx \omega_C^{\xx 2}\xx R_C\to 0.$$
Passing to the long exact sequence we have
$$0\to H^0\left(\bigwedge^i M_{\omega}\xx \omega_C\xx R_C\right)\to \bigwedge^i H^0(\omega_C)\xx H^0(\omega_C\xx R_C)\xrightarrow{\Lambda} H^0\left(\bigwedge^{i-1}M_{\omega}\xx \omega_C^{\xx2}\xx R_C\right)$$

\begin{prop}\label{prop_vir}
The two vector spaces 
$$\bigwedge^i H^0(\omega_C)\xx H^0(\omega_C\xx R_C)\ \m{ and }\ H^0\left(\bigwedge^{i-1}M_{\omega}\xx \omega_C^{\xx2}\xx R_C\right)$$
have the same dimension. As a consequence
any point $[C,\phi]$ of $\R_{g,S_3}^{S_3}$ is in $\mc U_g$
if and only if the associated $\Lambda$ morphism is not an isomorphism.
\end{prop}
\proof We start by proving that
\begin{equation}\label{eq_inizio}
\bigwedge^i H^0(\omega_C)\xx H^1(\omega_C\xx R_C)=H^1\left(\bigwedge^{i-1}M_{\omega}\xx \omega_C^{\xx 2}\xx R_C\right)=0.
\end{equation}
We observe that this two terms fit in the long exact sequence
$$\cdots\to \bigwedge^iH^0(\omega_C)\xx H^1(\omega_C\xx R_C)\to H^1\left(\bigwedge^{i-1}M_{\omega}\xx \omega_C^{\xx 2}\xx R_C\right)\to H^2\left(\bigwedge^iM_{\omega}\xx \omega_C\xx R_C\right)\to \cdots$$
and because the last term is $0$, it suffices to prove that $H^1(\omega_C\xx R_C)=H^0(R^{\vee})=0$.

Consider the group quotient $S_3\slash N=\mmu_2$, it
induces a twisted $\mmu_2$-cover $C\xrightarrow{\phi}BS_3\to B\mmu_2$.
This is equivalent to a principal $\mmu_2$-bundle $\varphi\colon F\to C$.
As $R_C$ is defined via the representation $\rho\colon S_3\to \GL(R)\cong\GL(2;\C)$,
its pullback on $F$ is given by the restriction $\rho|_N$ which is a representation of $N\cong \mmu_3$,
which decomposes as direct sum of two irreducible $\mmu_3$ representations of rank $1$.
We observe in particular that $\varphi^*R_C^{\vee}=\eta\oplus\eta^{\xx2}$ where $\eta$ is a non-trivial
line bundle such that $\eta^{\xx 3}\cong \Oo_F$. Therefore $H^0(\eta)=H^0(\eta^{\xx2})=0$ and
\cor{a fortiori} $H^0(R_C^{\vee})=0$.\newline

Since the first cohomology group is trivial, we have that the
Euler characteristic of both vector bundles coincides with the
dimension of their spaces of global sections. In particular, if
$\mc E$ is one of these vector bundles, $H^0(\mc E)=\chi(\mc E)=\deg(\mc E)+\rk(\mc E)(1-g)$.

\begin{itemize}
\item If $\mc E=\bigwedge^iH^0(\omega_C)\xx H^0(\omega_C\xx R_C)$, then
$$\dim (\mc E)=\dim\left(\bigwedge^iH^0(\omega_C)\right)\cdot \dim \left(H^0(\omega_C\xx R_C)\right).$$
The first term is simply $\binom{g}{i}$ because $h^0(\omega_C)=g$. For the second term we have 
$h^0(\omega_C\xx R_C)=\chi(\omega_C\xx R_C)=\deg(\omega_C\xx\R_C)+2-2g=2g-2$.
Therefore
$$\dim\left(\bigwedge^iH^0(\omega_C)\xx H^0(\omega_C\xx R_C)\right)=4i\binom{2i+1}{i}.$$

\item If $\mc E= H^0(\bigwedge^{i-1}M_{\omega}\xx \omega_C^{\xx 2}\xx R_C)$, then
$$\deg(\mc E)=2\cdot \deg\left(\bigwedge^{i-1}M_{\omega}\right)+\binom{g-1}{i-1}\cdot\deg(\omega_C^{\xx2}\xx R_C)=4(3i+1)\binom{2i}{i-1}.$$
Knowing that $\rk (\mc E)=2\cdot\binom{g-1}{i-1}$, we have
$$\dim(\mc E)=4\binom{2i}{i-1}(2i+1)=4i\binom{2i+1}{i},$$
and this completes the proof.
\end{itemize}\fine

\begin{prop}\label{propeff}
Consider a general hyperelliptic curve $(C,p)$ of genus $g\geq 2$ with a Weierstrass point.
Then there exists a twisted $S_3$-cover $\phi\colon C\to BG$ such that, if $R_C$ is the 
vector bundle of rank $2$ associated to the $R$ representation via $\phi$, then
$H^0(C,R_C((g-1)p))=0$.
\end{prop}

This proves the effectiveness of $\mc U_g$ by following the approach of \cite[Theorem 2.3]{cefs13}:
via the proposition we prove that $[C,\phi]$ is a point outside $\mc U_g$ in the case of genus $g=2i+1$, and therefore
that $\mc U_g$ is effective. Indeed, $M_{\omega_C}=\Oo_C(-2p)^{\oplus g-1}$ and therefore 
$$\bigwedge^{i}M_{\omega}\xx \omega_C\xx R_C={\left(R_C((g-1)p)\right)}^{\oplus\binom{g-1}{i}}.$$
This means that $[C,\phi]\notin \mc U_g$ if and only if $H^0(R_C((g-1)p))=0$.

\proof We recall that $R_C$ is obtained by pulling back via $\phi$ the irreducible representation $\rho\colon S_3\to\GL(R)\cong\GL(2;\C)$.
Consider the normal subgroup $N\subset S_3$, and the group quotient $S_3\slash N=\mmu_2$. This induces 
a twisted $\mmu_2$-cover $C\xrightarrow{\phi}BS_3\to B\mmu_2$. As $C$ is a scheme theoretic curve, this is simply a principal $\mmu_2$-bundle
$\pi\colon F\to C$. We use the notation $W:=R_C((g-1)p)$ and we want to show that
$\pi^*W$ decomposes as direct sum of two line bundles. By construction the image of $F\to C\to BS_3$ is $BN$, that is
$\pi^*W$ is induced by a representation of $N\cong \mmu_3$. The group $\mmu_3$ has $3$ irreducible
representations: the trivial one, the identity $\eta\colon \mmu_3\hookrightarrow \C^*$ and~$\eta^{\xx 2}$.
The restriction $\rho|_N$ decomposes as $\eta\oplus\eta^{\xx2}$. With a little abuse of notation we denote by
$\eta$ and $\eta^{\xx2}$ the order $3$ line bundles induced by these representations. Moreover,
we call $A_1:=\eta((g-1)p)$ and $A_2:=\eta^{\xx 2}((g-1)p)$.

By construction we have the decomposition $\pi^*W=\pi^*A_1\oplus\pi^*A_2$. We observe that
$H^0(W)=H^0(\pi^*W)_+$, that is the space of sections which are invariant by the
natural involution of $F$. If $H^0(\pi^*W)_+\neq 0$, then one between $H^0(\pi^*A_1)_+$ and $H^0(\pi^*A_2)_+$
is non-empty, but $H^0(\pi^*A_j)_+=H^0(A_j)$ and by \cite[Theorem 2.3]{cefs13}
there exists $(C,\phi)$ such that both $H^0(A_j)$ are empty. 

To prove this last point, observe that by \cite[Theorem 2.3]{cefs13}
there exists a point a point $[C,\nu]$ outside the divisor $\U_{g,3}\subset\overline\R_{g,\mmu_3}^{\mmu_3}$, where $\nu$ is a non-trivial
third root of $\Oo_C$, or equivalently a non-trivial principal $\mmu_3$-bundle on $C$. As the locus of hyperelliptic
curves with a Weiestrass point and a third root, is not included in $\U_{g,3}$, then by dimensional considerations
there exists an hyperelliptic curve $C$ such that $[C,\nu]\notin\U_{g,3}$ for any $\nu$. And this completes the proof.\fine

We use the result of Proposition \ref{prop_vir} to evaluate the class
of $\mc U_g$ in the Chow ring, and also the class of its closure
on the space $\widetilde \R_{g,S_3}^{S_3}$.
Consider the universal family $u\colon\cc_{g,S_3}\to\overline\R_{g,S_3}^{S_3}$
and the universal rank $2$ vector bundle $R_{\cc}$ associated to the representation $R$.
Introduce the vector bundle $\mfk M_u$ on $\overline\R_{g,S_3}^{S_3}$ defined by the short exact sequence
$$0\to \mfk M_u\to u^*(u_*\omega_u)\to\omega_u\to 0.$$
The geometric fiber of $\mfk M_u$ over any point $[\ssC,\phi]$ of $\R_{g,S_3}^{S_3}$ is
the previously defined vector bundle $M_{\omega}$.

\begin{defin}
We introduce a series of sheaves on $\overline\R_{g,S_3}^{S_3}$,
$$\E_{j,b}:=u_*\left(\bigwedge^j\mfk M_u\xx\omega_u^{\xx b}\xx R_{\cc}\right),$$
with $j\geq 0$ and $b\geq 1$ integers.
\end{defin}

\begin{rmk}
To prove that these sheaves are locally free on $\widetilde\R_{g,S_3}^{S_3}$,
by Grauert's Theorem 
it suffices to prove that $h^1(M_{\omega}\xx \omega_{\ssC}^{\xx b}\xx R_{\ssC})=0$
for every twisted $S_3$-cover $(\ssC,\phi)$ in $\widetilde\R_{g,S_3}^{S_3}$.
As showed in the proof of Proposition \ref{prop_vir}, this reduces
to prove $h^0((\omega_{\ssC}^{\vee})^{\xx (b-1)}\xx R_{\ssC}^{\vee})=0$ for all
$b\geq1$. In the proof of Proposition \ref{prop_vir} we showed this equality
for $b=1$ and a scheme theoretic curve $\ssC=C$, but the same
proof works for every $b\geq 1$ and for every twisted curve $\ssC$
such that its coarse space is irreducible.
\end{rmk}

By Proposition \ref{prop_vir}, on $\R_{g,S_3}^{S_3}$ the locus
$\mc U_g$ is the degeneration locus of a morphism
between the vector bundles
$\E_{i-1,2}$ and $\bigwedge^i\E\xx\E_{0,1}$,
where $\E$ is the Hodge bundle $u_*\omega_u$.
In the following we will evaluate this degeneration locus
on the space $\widetilde\R_{g,S_3}^{S_3}$.
 We denote by $\overline\U_g$ the closure of the locus $\U_g$
on the space $\widetilde \R_{g,S_3}^{S_3}$.

\begin{lemma}\label{l1}
Given an odd genus $g=2i+1$, the class of $\overline\U_g$ in $\Pic_\Q(\widetilde \R_{g,S_3}^{S_3})$ is
$$[\overline\U_g]=c_1\left(\Ho\left(\E_{i-1,2},\bigwedge^i\E\xx\E_{0,1}\right)\right)=
\rk(\E_{i-1,2})\cdot\left(\sum_{b=0}^i(-1)^{b+1}c_1\left(\bigwedge^{i-b}\E\xx\E_{0,b+1}\right)\right).$$
\end{lemma}
\proof Given two vector bundles $\mc A$ and $\mc B$ over $\widetilde \R_{g,S_3}^{S_3}$,
the morphism vector bundle $\Ho(\mc A,\mc B)$ is isomorphic to
$\mc B\xx \mc A^\vee$ and therefore
$c_1(\Ho(\mc A,\mc B))=\rk(\mc A)c_1(\mc B)-\rk(\mc B)c_1(\mc A)$.
To conclude we observe that by the definition of the vector bundles $\E_{j,b}$,
they fit in the short exact sequences
$$0\to \E_{j,b+1-j}\to\bigwedge^j \E\xx\E_{0,b+1-j}\to \E_{j-1,b+2-j}\to 0,$$
for all $j\geq 0$ and $b\geq j$.\fine

\begin{lemma}\label{l2}
The first Chern class of $\E_{0,b}$ is
$$c_1(\E_{0,b})=2 \lambda +2 \binom{b}{2}\kappa_1-\frac{1}{4}\delta_0^T-\frac{2}{3}\delta_0^N\in\Pic_\Q(\widetilde \R_{g,S_3}^{S_3}).$$
\end{lemma}
\proof This is a direct application of Proposition \ref{prop_chern} in
the evaluation of $\ch_1(Ru_*W_{\mc C})$.
By \cite{chio08.2} we have
$$\lambda=\ch_1(u_*\omega_u)=\frac{B_2}{2}\kappa_1+\frac{1}{2}\cdot\sum_{\substack{0\leq i\leq g-1,\\ 
\lq h\rq \in\lq G\rq } } r(h)\cdot(j_{i,\lq h\rq })_*\left(\frac{B_2}{2}\cdot\sum_{a+a'=d-1}\psi^a(-\psi')^{a'}\right).$$
To complete the proof
we only need the eigenvalues decomposition of $R((12))$ and $R((123))$,
where $R\colon S_3\to \GL(\C,2)$ is the irreducible $S_3$ representation
of dimension $2$.\fine

With these lemmata we can develop the calculations to evaluate
$[\overline\U_g]$ in terms of the Hodge class and the boundary classes.

\begin{prop}
In the Picard group $\Pic_\Q(\widetilde\R_{g,S_3}^{S_3})$ we have,
$$[\overline\U_g]=\rk(\E_{i-1,2})\cdot2\cdot\binom{2i-2}{i-1}\left(\frac{2(3i+1)}{i}\lambda-\delta_0'-\left(\frac{6i+1}{4i}\right)\delta_0^T-\left(\frac{5i+2}{3i}\right)\delta_0^N\right).$$
\end{prop}
\proof From the result of Lemma \ref{l2}, we
have that
$$c_1\left(\bigwedge^{i-b}\E\xx\E_{0,b+1}\right)=\rk(\bigwedge^{i-b}\E)\cdot c_1(\E_{0,b+1})+\rk(\E_{0,b+1})\cdot c_1(\bigwedge^{i-b}\E)=$$
$$=\binom{g}{i-b}\cdot\left(2\lambda+2\binom{b+1}{2}\kappa_1-\frac{1}{4}\delta_0^T-\frac{2}{3}\delta_0^N\right)+2(2b+1)(g-1)\cdot \binom{g-1}{i-b-1}\lambda.$$
This, thanks to Lemma \ref{l1}, allows
to conclude the evaluation.\fine

We are ready to introduce the divisor $\mm$ of Equation (\ref{eq_K})
over $\overline\R_{g,S_3}^{S_3}$, with $g=2i+1$ odd genus.
Harris and Mumford introduced in \cite{harmum82}
the following divisor.
\begin{defin}\label{def_mmm}
If $W^r_d(C)$ is the set of complete linear series over $C$
of degree $d$ and dimension at least $r$, it defines
the locus
$$\M_{g,i+1}^i:=\{[C]\in\M_g\ |\ W_{i+1}^1(C)\neq \varnothing\}\subset\M_g.$$
\end{defin}

In the same paper they proved
$$\M_{g,i+1}^i=c'\cdot\left(\frac{6(i+2)}{i+1}\lambda-\delta_0\right)\in\Pic_\Q(\widetilde\M_g),$$
where $\widetilde\M_g\subset\overline\M_g$ is the locus of irreducible stable curves,
and $c'$ is a positive coefficient.
Then, 
$$[ \pi^*\M_{g,i+1}^1]=c'\cdot\left(\frac{6(i+2)}{i+1}\lambda-\delta_0'-2\delta_0^T-3\delta_0^N\right)\in\Pic_Q(\widetilde \R_{g,S_3}^{S_3}),$$
where $\pi$ is the natural projection $\widetilde\R_{g,S_3}^{S_3}\to\widetilde \M_g$.\newline

Summarizing,
\begin{equation}\label{eq_kfin}
[K_{\widetilde \R_{g,S_3}}]=\alpha\cdot[\overline\U_g]+\beta\cdot[\pi^*\M_{g,i+1}^1]+E+\gamma\cdot\lambda,
\end{equation}
for every odd genus $g=2i+1$, with $\alpha$ and $\beta$ positive coefficients and $E$ an effective
sum of boundary divisors.

\begin{prop}\label{prop_fin}
In Equation (\ref{eq_kfin}) the $\gamma$ coefficient can be chosen strictly positive
for any~$i>5$.
\end{prop}
\proof By scaling appropriately every coefficient, the equation is equivalent
to choosing a real number $s\in [0,1]$ such that
$$s\cdot\left(\frac{2(3i+1)}{i}\lambda-\delta_0'-\left(\frac{6i+1}{4i}\right)\delta_0^T-\left(\frac{5i+2}{3i}\right)\delta_0^N\right)+$$
$$+(1-s)\cdot\left(\frac{6(i+2)}{i+1}\lambda-\delta_0'-2\delta_0^T-3\delta_0^N\right) + E +\gamma \cdot \lambda=$$
$$=\frac{13}{2}\lambda-\delta_0'-\frac{3}{2}\delta_0^T-2\delta_0^N.$$
For $E$ to be an effective divisor we must have
$s\cdot\left(\frac{6i+1}{4i}\right)+(1-s)\cdot 2\geq \frac{3}{2}$
and
$s\cdot\left(\frac{5i+2}{3i}\right)+(1-s)\cdot 3\geq 2$,
and therefore $s\leq \frac{3i}{4i-2}$ is a necessary and sufficient
condition for $E$ to be an effective boundary divisor.

To complete the proof we evaluate the $\gamma$ coefficient,
$$s\cdot \left(\frac{6i+2}{i}\right)+(1-s)\cdot \left(\frac{6i+12}{i+1}\right)+\gamma=\frac{13}{2}.$$
After calculations this gives
$\gamma=\frac{i-11}{2(i+1)}+s\cdot\frac{4i-2}{i(i+1)}$,
which means a maximal possible value of 
$$\gamma=\frac{i-5}{2(i+1)},$$
which is positive if and only if $i>5$.\fine

After this proposition, 
the canonical divisor over $\widetilde\R_{g,S_3}^{S_3}$
is big for every odd genus $g\geq 13$. Then, thanks to Proposition \ref{prop_big}, the moduli space
 $\overline\R_{g,S_3}^{S_3}$ is of general type for every odd genus $g\geq 13$.

\bibliography{mybiblio}{}
\bibliographystyle{plain}

\end{document}